\documentclass[journal]{IEEEtran}

\usepackage{graphicx}

\usepackage{amsmath}
\DeclareMathOperator{\Tr}{Tr}
\usepackage{amssymb}
\usepackage{amsthm}
\usepackage{url}

\usepackage[utf8]{inputenc}
\usepackage[T1]{fontenc}
\usepackage{booktabs}
\usepackage{multirow}
\usepackage{makecell}
\usepackage{xcolor}
\usepackage[binary-units]{siunitx}
\usepackage{interval}

\newtheorem{lem}{Lemma}

\newtheorem{thm}{Theorem}
\newtheorem{rem}{Remark}
\newtheorem{ass}{Assumption}
\newtheorem{cor}{Corollary}
\begin{document}
%
\title{Efficient Representation and Approximation of Model Predictive Control Laws via Deep Learning}
%
%
%

\author{Benjamin~Karg,
        Sergio~Lucia,~\IEEEmembership{Member,~IEEE,}
\thanks{B. Karg and S. Lucia are with the Chair of Internet of Things for Smart Buildings, TU Berlin, and Einstein Center Digital Future, Einsteinufer 17, 10587 Berlin, Germany, e-mail: benjamin.karg@tu-berlin.de, sergio.lucia@tu-berlin.de.}}

\markboth{IEEE Transactions on Cybernetics}%
{Shell \MakeLowercase{\textit{et al.}}: Bare Demo of IEEEtran.cls for IEEE Journals}
%



\onecolumn
© 2020 IEEE.  Personal use of this material is permitted.  Permission from IEEE must be obtained for all other uses, in any current or future media, including reprinting/republishing this material for advertising or promotional purposes, creating new collective works, for resale or redistribution to servers or lists, or reuse of any copyrighted component of this work in other works.
\twocolumn
\maketitle

\begin{abstract}

We show that artificial neural networks with rectifier units as activation functions can exactly represent the piecewise affine function that results from the formulation of model predictive control of linear time-invariant systems.
The choice of deep neural networks  is particularly interesting as they can represent exponentially many more affine regions compared to networks with only one hidden layer.
We provide theoretical bounds on the minimum number of hidden layers and neurons per layer that a neural network should have to exactly represent a given model predictive control law.

The proposed approach has a strong potential as an approximation method of predictive control laws, leading to better approximation quality and significantly smaller memory requirements than previous approaches, as we illustrate via simulation examples.
We also suggest different alternatives to correct or quantify the approximation error.
Since the online evaluation of neural networks is extremely simple, the approximated controllers can be deployed on low-power embedded devices with small storage capacity, enabling the implementation of advanced decision-making strategies for complex cyber-physical systems with limited computing capabilities.

\end{abstract}


\IEEEpeerreviewmaketitle

\section{Introduction}
Model predictive control (MPC) is a popular control strategy that computes control inputs by solving a numerical optimization problem. A mathematical model is used to predict the future behavior of the system and an optimal sequence of  control inputs is computed by solving an optimization problem that minimizes a given objective function subject to constraints. The main reasons for its success are the possibility of handling systematically multiple-input multiple-output systems, nonlinearities as well as constraints.
The main challenge of MPC is that it requires the solution of an optimization problem at each sampling time of the controller. For this reason, traditional applications included those related to slow systems such as chemical processes~\cite{qin2003}, \cite{rawlings2009}. 

During the past two decades, a large research effort has been devoted to broaden the range of MPC applications, leading to various specific MPC algorithms reaching from event-triggered robust control~\cite{liu2018aperiodic} to hierarchical distributed systems ~\cite{hans2019hierarchical}.
To enable the application of MPC strategies to complex cyber-physical systems \cite{raman2014}, extending the application of MPC to fast embedded systems with limited computing capabilities is an important challenge.
Two different approaches have been followed to achieve this goal.
The first approach included the development of fast solvers and tailored implementations  \cite{mattingley2012cvxgen} that can solve the required optimization problems in real time for fast systems. Different variations of the Nesterov's fast gradient method (see e.g. \cite{richter2012} \cite{giselsson2013}, \cite{koegel2011}) and of the alternating directions method of multipliers (ADMM)  \cite{boyd2011distributed} have been very successful for embedded optimization and model predictive control.  Different versions of these algorithms have been used to obtain MPC implementations on low-cost microcontrollers  \cite{zometa2013}, \cite{lucia2016_ARC} or high-performance FPGAs \cite{jerez2014}, \cite{lucia2018_TII}.

The second approach to extend the application of MPC to fast and embedded systems is usually called explicit MPC. The MPC problem for linear time invariant systems is a parametric quadratic program whose solution is a piecewise affine function defined on polytopes and only depends on the current state of the system \cite{bemporad2002}. Explicit MPC exploits this idea by precomputing and storing the piecewise affine function that completely defines the MPC feedback law. The online evaluation of the explicit MPC law reduces to finding the polytopic region in which the system is currently located and applying the corresponding affine law.
The main drawback of explicit MPC is that the number of regions on which the control law is defined grows exponentially with the prediction horizon and the number of constraints.
The inherent growth of the memory footprint and of the complexity of the point location problem limits the application to small systems and small prediction horizons, especially in the case of embedded systems with limited storage capabilities and computational power.

To counteract the massive memory requirements, some approaches try to simplify the representation of the control law by eliminating redundant regions \cite{geyer2008} or by using different number representations \cite{ingole2017}. Other approaches try to approximate the exact explicit MPC solution to further reduce the memory requirements of the approach (see a review in \cite{alessio2009}). Approximate explicit MPC schemes include the use of simplicial partitions \cite{Bemporad2011}, neural networks \cite{parisini1995}, radial basis functions \cite{cseko2015}, lattice representation \cite{wen2009analytical} or using a smaller number of regions to describe the MPC law \cite{holaza2013}.

To reduce the complexity of the point location problem, binary search trees (BST) \cite{tondel2003evaluation} introduce a tree structure where the nodes represent unique hyperplanes. At each node it is checked on which side of the hyperplane the state is until a leaf node is reached. At the leaf node, a unique feedback law is identified and evaluated.
This method renders the online computation time logarithmic in the number of regions, but precomputation times can be prohibitive or intractable for larger problems \cite{bayat2012flexible}.
Modifications of BST include approximations of the exact solution via hypercubic regions \cite{johansen2003approximate}, truncated BSTs combined with direct search  to restrict the depth of BSTs \cite{bayat2011combining}, computation of arbitrary hyperplanes which balance the tree and minimize its depth \cite{fuchs2010optimized} and merging an BST with the lattice representation \cite{bayat2012flexible}.


Motivated by new advances on the theoretical description of the representation capabilities of deep neural networks \cite{silver2016}, \cite{safran2017}, the goal of this work is to provide a scheme for an approximated explicit MPC controller with significantly lower memory and computational requirements compared to existing approaches.
Deep neural networks (with several hidden layers) can represent exponentially many more linear regions than shallow networks (with only one hidden layer).
This attribute has been exploited in recent works for complex control tasks such as mixed-integer MPC \cite{karg2018} and robust nonlinear MPC \cite{lucia2018deep}.
While many control approaches have used neural networks to capture unknown or nonlinear dynamics of the model~\cite{liu2018neural}, \cite{ding2018neural}, in this work neural networks are used to directly approximate the optimal control law.
This work differs from \cite{chen2018approximating}, which is based on similar ideas, by presenting bounds on the necessary size of a deep network to represent an explicit MPC law exactly, as well as by presenting a comprehensive comparison with other state-of-the-art approximate explicit MPC methods.
Also, statistical verification techniques are provided that can be evaluated with high-fidelity models, even if the controller was designed with simpler models.
The main contributions of the paper are:
\begin{itemize}
	\item the derivation of explicit bounds for the required size (width and depth) that a deep network should have to exactly represent a given explicit MPC solution.
	\item the presentation of an approach to approximate explicit MPC based on deep learning which achieves better accuracy with less memory requirements when compared to other approximation techniques.
	\item statistical verification techniques to assess the validity of the obtained approximate controllers.
	\item an embedded implementation of the resulting controllers.
\end{itemize}

The remainder of the paper is organized as follows. Section \ref{sec:background} introduces background information about model predictive control and neural networks. Section \ref{sec:deep} presents explicit bounds for a deep network to be able to represent exactly an MPC law and serves as a motivation for the approximation of explicit MPC laws presented in Section \ref{sec:approximate}, where different techniques to deal with the approximation error are also presented. Section \ref{sec:results} illustrates the potential of the approach with two simulation examples and the paper is concluded in Section \ref{sec:conclusions}.
\section{Background and Motivation}\label{sec:background}
\subsection{Notation}
We denote by $\mathbb{R}$, $\mathbb{R}^{n}$ and $\mathbb{R}^{n \times m}$ the real numbers, $n$-dimensional real vectors and $n \times m$ dimensional real matrices, respectively.
The interior of a set  is denoted by $\text{int}(\cdot)$, its cardinality by $|\cdot|$ and $\left \lfloor \cdot \right \rfloor$ denotes the \emph{floor} operation, i.e. the rounding to the nearest lower integer.
The probability of an event is denoted by $P(\cdot)$ and the composition of two functions $f$ and $g$ by $g\circ f(\cdot) = g(f(\cdot))$.

\subsection{Explicit MPC}
Model predictive control (MPC) is an optimal control scheme that uses a system model to predict the future evolution of a system. 
We consider discrete linear time-invariant (LTI) systems:
\begin{align}\label{eq:LTI}
x_{k+1} = Ax_k + Bu_k,
\end{align}
where $x \in \mathbb{R}^{n_x}$ is the state vector, $u \in \mathbb{R}^{n_u}$ is the control input, $A \in \mathbb{R}^{n_x \times n_x}$ is the system matrix, $B \in \mathbb{R}^{n_x \times n_u}$ is the input matrix and the pair $(A,B)$ is controllable.

Using a standard quadratic cost function, the following constrained finite time optimal control problem with a horizon of $N$ steps should be solved at each sampling time to obtain the MPC feedback law:

\begin{subequations}\label{eq:FTOCP}
	\begin{align}
		&\underset{\tilde u}{\text{minimize}} && x_N^T P x_N + \sum_{k=0}^{N-1}{x_k^T Q x_k + u_k^T R u_k} \span \span \span \\
		&\text{subject to} &&x_{k+1} = Ax_k + Bu_k, \\
		&            &&C_x x_k \leq c_x,\, C_f x_N \leq c_f, \label{eq:state_terminal_cons}\\
		&            &&C_u u_k \leq c_u, \label{eq:input_cons}\\
		&			 &&x_0 = x_{\text{init}}, \\
		& && \forall \,\, k = 0,\dots, N-1,
	\end{align}
\end{subequations}
where $\tilde u = [u_0, \dots, u_{N-1}]^T$ is a vector that contains the sequence of control inputs and $P \in \mathbb{R}^{n_x \times n_x}$, $Q \in \mathbb{R}^{n_x \times n_x}$ and $R \in \mathbb{R}^{n_u \times n_u}$ are the weighting matrices. The weighting matrices are chosen such that $P \succeq 0$ and $Q \succeq 0$ are positive semidefinite, and $R \succ 0$ is positive definite. The state, terminal and input constraints are bounded polytopic sets $\mathcal X$, $\mathcal X_f$ and $\mathcal U$ defined by the matrices $C_x\in \mathbb R^{n_{\text{cx}}\times n_x}$, $C_f\in \mathbb R^{n_{\text{cf}}\times n_x}$, $C_u\in \mathbb R^{n_{\text{cu}}\times n_u}$ and the vectors $c_x\in \mathbb R^{n_{\text{cx}}}$, $c_f\in \mathbb R^{n_{\text{cf}}}$, $c_u\in \mathbb R^{n_{\text{cu}}}$.
The terminal cost defined by $P$ as well as the terminal set $\mathcal X_f$ are usually chosen in such a way that stability of the closed-loop system and recursive feasibility of the optimization problem are guaranteed~\cite{mayne2000}. The set of initial states $x_{\text{init}}$ for which~\eqref{eq:FTOCP} has a feasible solution depending on the prediction horizon $N$ is called feasibility region and is denoted by $\mathcal X_N$.

The optimization problem \eqref{eq:FTOCP} can be reformulated as a multi-parametric problem \cite{bemporad2002} that only depends on the current system state $x_{\text{init}}$:
\begin{subequations}\label{eq:condensed}
	\begin{align}
		&\underset{\tilde u}{\text{minimize}}&& \tilde u^T F \tilde u + x_{\text{init}}^T G \tilde u + x_{\text{init}}^T H x_{\text{init}} \span \\
		&\text{subject to} && C_c \tilde u \leq T x_{\text{init}} + c_c,
	\end{align}
\end{subequations}
where $F \in \mathbb{R}^{N n_u \times N n_u}$, $G \in \mathbb{R}^{n_x \times N n_u}$, $H \in \mathbb{R}^{n_x \times n_x}$, $C_c \in \mathbb{R}^{N n_{\text{ineq}} \times N n_u}$, $T \in \mathbb{R}^{N n_{\text{ineq}} \times n_x}$, $c_c \in \mathbb{R}^{N n_{\text{ineq}}}$ and $n_{\text{ineq}}$ is the total number of inequalities in~\eqref{eq:FTOCP}.

The solution of the multi-parametric quadratic programming problem~\eqref{eq:condensed} is a piecewise affine (PWA) function of the form \cite{bemporad2002}:
\begin{align}
\mathcal{K}(x_{\text{init}}) = 
	\begin{cases}
		K_1 x_{\text{init}}+g_1 & \text{if} \quad x_{\text{init}} \in \mathcal{R}_1, \\
		& \vdots \\
		K_{n_{\text{r}}} x_{\text{init}}+g_{n_{\text{r}}} & \text{if} \quad x_{\text{init}} \in \mathcal{R}_{n_{\text{r}}},
	\end{cases}
\label{eq:exp_mpc}
\end{align}
with $n_{\text{r}}$ regions, $K_i \in \mathbb{R}^{Nn_u \times n_x}$ and $g_i \in \mathbb{R}^{Nn_u}$.
Each region $\mathcal R_i$ is described by a polyhedron
\begin{align}
\mathcal{R}_i = \{x \in \mathbb{R}^{n_x} \mid Z_i x \leq z_i\} \quad \forall i = 1,\dots,n_r,
\end{align}
where $Z_i \in \mathbb{R}^{c_i \times n_x}$, $z_i \in \mathbb{R}^{c_i}$ describe the $c_i$ halfspaces $a_{i,j} x_{\text{init}} \leq b_{i,j}$ of the $i$-th region with $j = 1,\dots,c_i$, $a_{i,j} \in \mathbb{R}^{1 \times n_x}$ and $b_{i,j} \in \mathbb{R}$.
The formulation \eqref{eq:exp_mpc} is defined on the bounded polytopic partition $\mathcal{R}_\Omega = \cup_{i=1}^{n_{\text{r}}} \mathcal{R}_i$ with $\text{int}(\mathcal{R}_i) \cap \text{int}(\mathcal{R}_j) = \emptyset$ for all $i \neq j$.

Most hyperplanes are shared by neighbouring regions and the feedback law can be identical for two or more regions.
Hence, the memory needed to store the explicit MPC controller~\eqref{eq:exp_mpc} can be approximated as
\begin{align}\label{eq:memory_empc}
\Gamma_{\mathcal{K}} = \alpha_{\text{bit}} \left( \left( n_{\text{h}} \left( n_x + 1 \right) \right) + n_{\text{f}} \left(n_x n_u + n_u \right) \right),
\end{align}
where $n_{\text{h}}$ is the number of unique hyperplanes, $n_{\text{f}}$ is the number of unique feedback laws and $\alpha_{\text{bit}}$  is the memory necessary to store a real number. Since for the actual implementation of the explicit MPC law only the input of the first time step is needed, only the  first $n_u$ rows of $K_j$ and $g_j$ for $j=1,\dots,n_{\text{f}}$ have to be stored which equals $n_x n_u+n_u$ numbers per unique feedback law.

One main drawback of the explicit MPC formulation is that the number of regions for an exact representation can grow exponentially with respect to the horizon and number of constraints~\cite{bemporad2002}, which leads to large memory requirements and might render the application of the method intractable.
For the simple example of the inverted pendulum on a cart, which is presented in detail in Section~\ref{sec:results}, with $n_x=4$ states, $n_u=1$ input and box constraints on the control input and the states, the explicit solution consists of $n_{\text{r}} = 91$ regions for $N=3$, $n_{\text{r}} = 191$ regions for $N=4$ and $n_{\text{r}} = 323$ regions for $N=5$. For $N=10$ as many as $n_{\text{r}} = 1638$ regions are obtained.

\subsection{Artificial Neural Networks}
This subsection shortly recaps the fundamental concepts of artificial neural networks.
A feed-forward neural network is defined as a sequence of layers of neurons which determines a function $\mathcal{N}:\mathbb{R}^{n_x} \rightarrow \mathbb{R}^{n_u}$ of the form
\begin{equation}\label{eq:neural_network}
\begin{split}
\mathcal{N}(x;\theta,M,L) = \qquad \qquad \qquad \qquad \qquad \qquad \qquad \qquad\\
\quad \bigg \{
\begin{array}{lll}
		f_{L+1} \circ g_L \circ f_L \circ \dots \circ g_1 \circ f_1(x) & \text{for} & L \geq 2, \\
		f_{L+1} \circ g_1 \circ f_1(x), & \text{for} & L = 1,
\end{array}
\end{split}
\end{equation}
where the input of the network is $x \in \mathbb{R}^{n_x}$ and the output of the network is $ u \in \mathbb{R}^{n_u}$.
$M$ is the number of neurons in each hidden layer and $L$ is the number of hidden layers.
If $L\geq2$, $\mathcal N$ is described as a \emph{deep} neural network and if $L=1$ as a \emph{shallow} neural network.
Each hidden layer consists of an affine function:
\begin{align}\label{eq:affine_function}
f_l(\xi_{l-1}) = W_l\xi_{l-1}+b_l,
\end{align}
where $\xi_{l-1} \in \mathbb{R}^M$ is the output of the previous layer with $\xi_0 = x$. The second element of the neural network is a nonlinear activation function $g_l$.
In this paper, exclusively rectifier linear units (ReLU) are considered as activation function, which compute the element-wise maximum between zero and the affine function of the current layer $l$:
\begin{align}\label{eq:ReLU}
g_{l}(f_{l}) = \max(0,f_{l}). 
\end{align}
The parameter $\theta = \{\theta_1, \dots, \theta_{L+1}\}$ contains all the weights and biases of the affine functions of each layer
\begin{align}
\theta_l = \{W_l,b_l\} \quad \forall l = 1, \dots, L+1,
\end{align}
where the weights are
\begin{align}
W_l \in
	\begin{cases}
		\mathbb{R}^{M \times n_x} & \text{if} \quad l=1, \\
		\mathbb{R}^{M \times M} & \text{if} \quad l=2,\dots,L, \\
		\mathbb{R}^{n_u \times M} & \text{if} \quad l=L+1,
	\end{cases}
\end{align}
and the biases are
\begin{align}
b_l \in
	\begin{cases}
		\mathbb{R}^M & \text{if} \quad l=1,\dots,L, \\
		\mathbb{R}^{n_u} & \text{if} \quad l = L+1.
	\end{cases}
\end{align}

\subsection{Motivation}

Implementing the explicit model predictive controller defined by \eqref{eq:exp_mpc} requires storing the unique sets of hyperplanes and feedback laws that define the regions and the affine controllers in each region~\eqref{eq:memory_empc}, which was shown before to be possibly prohibitive due to the exponential growth of the number of regions.

The main motivation of this work is to find an efficient representation and approximation of the MPC control law, which can significantly reduce the memory requirements for the representation of the exact controller as well for achieving a high-quality approximation that outperforms other approximate explicit MPC techniques.
 
The basic idea that motivates this work is described in the following Lemma.
\begin{lem}{\cite{montufar2014}}\label{lemma:expo}
Every neural network $\mathcal N(x;\theta, M, L)$ with input $x\in \mathbb R^{n_x}$ defined as in \eqref{eq:neural_network} with ReLUs as activation functions and $M \geq n_x$ represents a piecewise affine function. In addition, a lower bound on the maximal number of affine regions that the neural network represents is given by the following expression:
\begin{align*}
\left( \prod \limits_{l=1}^{L-1} \left \lfloor \frac{M}{n_x}\right \rfloor^{n_x} \right) \sum \limits_{j=0}^{n_x} \binom{L}{j}.
\end{align*}
\end{lem}

\begin{proof}[Proof of Lemma~\ref{lemma:expo}]
The neural network $\mathcal N(x;\theta, M, L)$ is a piecewise affine function because it only contains compositions of affine transformations with a piecewise affine function (ReLUs).
For the derivation of the maximal number of regions, see \cite{montufar2014}.
\end{proof}

Lemma~\ref{lemma:expo} gives clear insights about why deep networks, as often observed in practice, obtain better performance to approximate complex functions when compared to shallow networks. In particular, Lemma~\ref{lemma:expo} implies that the number of affine regions that a neural network can represent grows exponentially with the number of layers $L$ as long as the width of the network $M$ is not smaller than the number of inputs $n_x$. The bound of Lemma~\ref{lemma:expo} can be slightly improved if $M\geq 3n_x$ as recently shown in \cite{serra2018}.

\begin{figure}[t]
\begin{center}
\includegraphics[width=0.9\columnwidth]{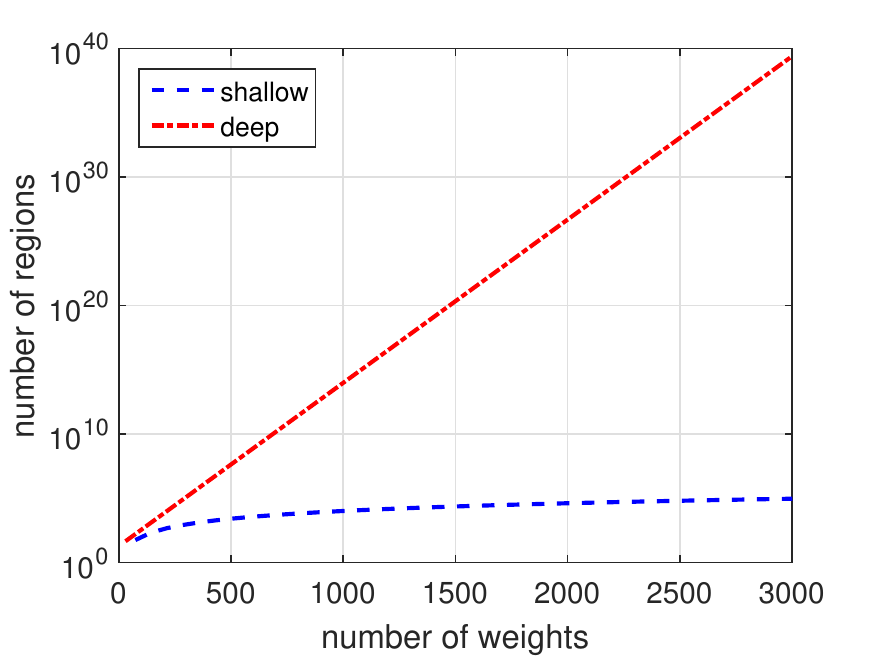} 
\caption{Number of regions with respect to the number of weights a neural network can represent. The parameters for this plot were chosen to $n_x = 2$, $M = 10$, $L = 1, \dots, 50$ and $n_u = 4$.}
\label{fig:regions_weights}
\end{center}
\end{figure}

At the same time, the number of parameters contained in $\theta$ that are necessary to fully describe the neural network $\mathcal N(x;\theta, M, L)$ are determined by the dimensions of the weights and biases at each layer. Assuming that storing each number requires $\alpha_{\text{bit}}$ bits, the total amount of memory necessary to store the neural network $\mathcal N(x;\theta, M, L)$ can be computed as:
\begin{align}\label{eq:mem_deep}
\Gamma_{\mathcal{N}} =\alpha_{\text{bit}}(&(n_x+1)M + (L-1)(M+1)M\notag \\& + (M+1)n_u).
\end{align}
Since $\Gamma_{\mathcal{N}}$ only grows linearly with respect to the number of layers $L$, deep ReLU networks can represent exponentially many more linear regions than shallow ones for a fixed amount of memory. This fact can be clearly seen in Fig.~\ref{fig:regions_weights}.

We believe that this observation, while somewhat obvious, is a very powerful result with important implications in control theory and constitutes the main motivation for this paper.
\section{Deep learning-based explicit MPC}\label{sec:deep}

This section shows how to design a deep neural network that can exactly represent an explicit MPC feedback law~\eqref{eq:exp_mpc} of the form $\mathcal{K}(x):[0,1]^{n_x} \rightarrow {\mathbb{R}^+}^{n_u}$ by only mapping the state $x$ to the first optimal control input $u_0^*$.
Considering only the control input of the first time step is sufficient, because after its application a new control input trajectory $\tilde{u}^*$ is computed in the MPC setting.
%
We make use of two lemmas to derive specific bounds for a deep neural network to be able to exactly represent any explicit MPC law $\mathcal{K}(x)$.
The following lemma from \cite{kripfganz1987} is used.
\begin{lem}\label{lemma:decomp}
Every scalar PWA function $f(x): \mathbb{R}^{n_x} \rightarrow \mathbb{R}$ can be written as the difference of two convex PWA functions
\begin{align}
f(x) = \gamma(x) - \eta(x),
\end{align}
where $\gamma(x): \mathbb{R}^{n_x} \rightarrow \mathbb{R}$ has $r_\gamma$ regions and $\eta(x): \mathbb{R}^{n_x} \rightarrow \mathbb{R}$ has $r_\eta$ regions.
\end{lem}
\begin{proof}[Proof of Lemma~\ref{lemma:decomp}]
See \cite{kripfganz1987} or \cite{hempel2013}.
\end{proof}
The following Lemma, recently presented in \cite{hanin2017}, gives specific bounds for the structure that a deep neural network should have to be able to exactly represent a convex piecewise affine function.
\begin{lem}\label{lemma:maximum}
A convex piecewise affine function $f:[0,1]^{n_x}\rightarrow \mathbb{R}^+$ defined as the point-wise maximum of $N$ affine functions:
\begin{align*}
f(x) = \max_{i=1,\dots, N}f_i(x),
\end{align*}
can be exactly represented by a deep ReLU network with width $M = n_x + 1$ and depth $N$.
\end{lem}
\begin{proof}[Proof of Lemma~\ref{lemma:maximum}]
See Theorem 2 from \cite{hanin2017}.
\end{proof}
One of the main contributions of this paper is given in the following theorem, which states that any explicit MPC law of the form~\eqref{eq:exp_mpc} with $\mathcal{K}(x):[0,1]^{n_x} \rightarrow {\mathbb{R}^+}^{n_u}$ can be represented by a deep ReLU neural network with a predetermined size.
\begin{thm}
\label{thm:main_theorem}
There always exist parameters $\theta_{\gamma,i}$ and $\theta_{\eta,i}$ for $2n_u$ deep ReLU neural networks with depth $r_{\gamma,i}$ and $r_{\eta,i}$ for $i = 1,\dots,n_u$ and width $M = n_x + 1$, such that the vector of neural networks defined by
\begin{align}
\label{eq:theorem}
\begin{split}
&\begin{bmatrix}
\mathcal{N}(x;\theta_{\gamma,1},M,r_{\gamma,1}) - \mathcal{N}(x;\theta_{\eta,1},M,r_{\eta,1}) \\
\vdots \\
\mathcal{N}(x
;\theta_{\gamma,n_u},M,r_{\gamma,n_u}) - \mathcal{N}(x;\theta_{\eta,n_u},M,r_{\eta,n_u}) \\
\end{bmatrix} 
\end{split}
\end{align}
can exactly represent an explicit MPC law $\mathcal{K}(x): [0,1]^{n_x} \rightarrow {\mathbb{R}^+}^{n_u}$.
\end{thm}

\begin{proof}[Proof of Theorem~\ref{thm:main_theorem}]
Every explicit MPC law $\mathcal{K}(x): [0,1]^{n_x} \rightarrow \mathbb{R}^{+^{n_u}}$ can be split into one explicit MPC law per output dimension:
\begin{align}
\quad \mathcal{K}_{i}(x): [0,1]^{n_x} \rightarrow \mathbb{R}^+ \quad \forall i = 1,\dots,n_u.
\end{align}
Applying Lemma~\ref{lemma:decomp} to all $n_u$ MPC laws, each one of them can be decomposed into two convex scalar PWA functions:
\begin{align}
\label{eq:proof_decomposed_nonvector}
\mathcal{K}_{i}(x) = \gamma_i(x) - \eta_i(x) \quad \forall i = 1,\dots,n_u,
\end{align}
where each $\gamma_i(x)$ and each $\eta_i(x)$ are composed of $r_{\gamma_i}$ and $r_{\gamma_i}$ affine regions.
The explicit MPC law $\mathcal{K}(x): [0,1]^{n_x} \rightarrow {\mathbb{R}^+}^{n_u}$ can thus be vectorized as
\begin{align}
\label{eq:proof_decomp}
\mathcal{K}(x) = 
\begin{bmatrix}
\gamma_1(x) - \eta_1(x) \\
\vdots \\
\gamma_{n_u}(x) - \eta_{n_u}(x) \\
\end{bmatrix}.
\end{align}
According to Lemma~\ref{lemma:maximum}, it is always possible to find parameters $\theta_{\gamma,i}$, $\theta_{\eta,i}$ for  deep ReLU networks with width $M=n_x+1$, depth not larger than $r_{\gamma,i}$, $r_{\eta,i}$  that can exactly represent the scalar convex functions $\gamma_i(x)$ and $\eta_i(x)$.
This holds because any convex affine function with $N$ regions can be described as the point-wise maximum of $N$ scalar affine functions.
This means that each component of the transformed explicit MPC law can be written as:
\begin{align}
\label{eq:proof_NN}
\begin{split}
\gamma_{i}(x) - \eta_{i}(x) &=\mathcal{N}(x;\theta_{\gamma,i},M,r_{\eta,i}) - \mathcal{N}(x;\theta_{\eta,i},M,r_{\eta,i}), \\
\end{split}
\end{align}
for all $i = 1,\dots,n_u$.
Substituting \eqref{eq:proof_NN} in \eqref{eq:proof_decomp} results in \eqref{eq:theorem}.
\end{proof}

Theorem~\ref{thm:main_theorem} requires that the piecewise affine function maps the unit hypercube to the space of positive real numbers.
Any explicit MPC law~\eqref{eq:FTOCP} can be written in this form, as long as the invertible affine transformations defined in Assumption~\ref{ass:transform} exist.
This result is formalized in the Corollary~\ref{cor:cor}.
\begin{ass}
\label{ass:transform}
There exist two invertible affine transformations $\mathcal{A}_x: \mathcal{R}_\Omega \rightarrow [0,1]^{n_x}$ and $\mathcal{A}_u: \mathcal{U} \rightarrow {\mathbb{R}^+}^{n_u}$ for an explicit control law~\eqref{eq:exp_mpc} $\mathcal{K}_{\text{orig}}: \mathcal{R}_\Omega \rightarrow \mathcal{U}$ such that
\begin{align} \label{eq:exp_mpc_trans}
\mathcal{K}_{\text{orig}}(x) = \mathcal{A}_u^{-1} \circ \mathcal{K}(\hat{x}),
\end{align}
where $\hat{x} = \mathcal{A}_x \circ x$. The affine transformations $\mathcal{A}_x$ and $\mathcal{A}_u$ always exist, when $\mathcal{R}_{\Omega}$ and $\mathcal{U}$ are compact sets, as it is standard in control applications.
\end{ass}
\begin{cor}\label{cor:cor}
If for a given explicit MPC solution~\eqref{eq:exp_mpc} $\mathcal{K}_{\text{orig}}: \mathcal{R}_\Omega \rightarrow \mathcal{U}$, there exist two invertible affine transformations such that Assumption~\ref{ass:transform} holds.
Theorem~\ref{thm:main_theorem} can be applied to the transformed MPC solution $\mathcal{K}(\hat{x}): [0,1]^{n_x}  \rightarrow {\mathbb{R}^+}^{n_u}$~\eqref{eq:exp_mpc_trans}.
Hence, such an explicit MPC solution of the form~\eqref{eq:exp_mpc} can be exactly represented by two invertible affine transformations and $2n_u$ deep ReLU networks with width and depth as defined in Theorem~\ref{thm:main_theorem}.
\end{cor}

The proof presented in \cite{hempel2013} for the decomposition of a PWA function into the difference of two PWA functions is constructive, which means that Theorem~\ref{thm:main_theorem} gives explicit bounds for the construction of neural networks that can exactly represent any explicit MPC of the form~\eqref{eq:exp_mpc} with $\mathcal{K}(x):[0,1]^{n_x} \rightarrow {\mathbb{R}^+}^{n_u}$, considering only the first step of the optimal control input sequence.

Another advantage of the proposed approach is that if the explicit MPC law is represented as a set of neural networks, its online application does not require determining the current region $\mathcal{R}_i$ and only needs the evaluation of the neural networks. This evaluation is a straightforward composition of affine functions and simple nonlinearities, which facilitates the implementation of the proposed controller on embedded systems.
 
We illustrate Theorem~\ref{thm:main_theorem} with a small example of an oscillator with the discrete system matrices

\begin{equation*}
\begin{aligned}
&A=
\begin{bmatrix}
0.5403   & 0.8415 \\ 
0.8415   & 0.5403 \\
\end{bmatrix},
&B=
\begin{bmatrix}
-0.4597 \\
0.8415 \\
\end{bmatrix}.
\end{aligned}
\end{equation*}

\begin{figure*}
  \includegraphics[width=.24\textwidth]{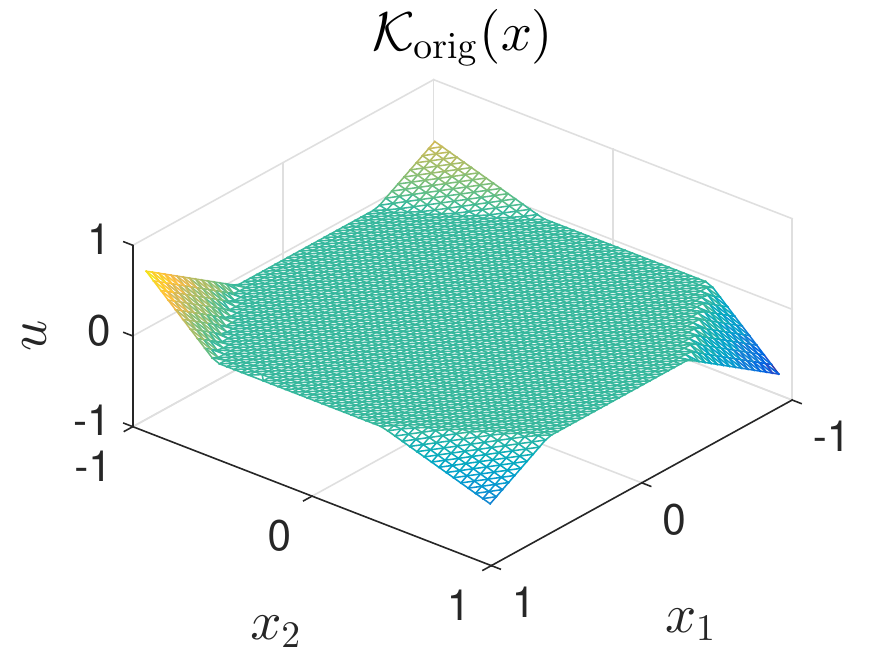}
  \includegraphics[width=.24\textwidth]{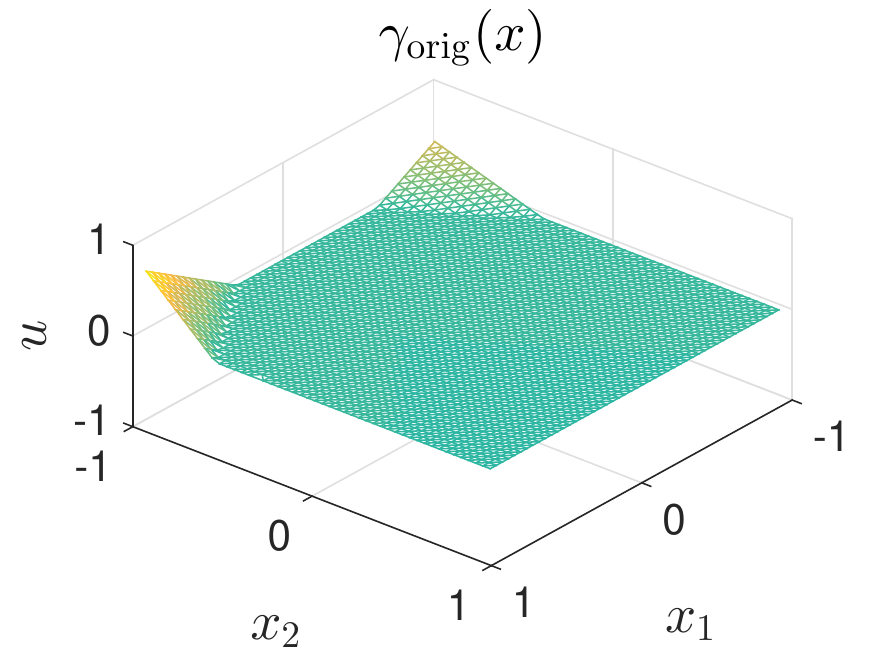}
  \includegraphics[width=.24\textwidth]{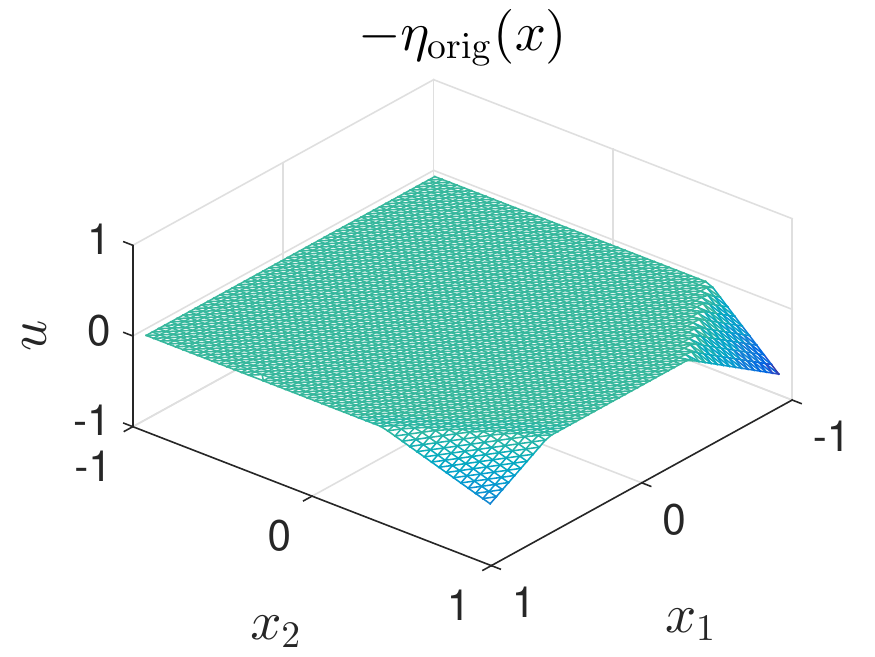}
  \includegraphics[width=.24\textwidth]{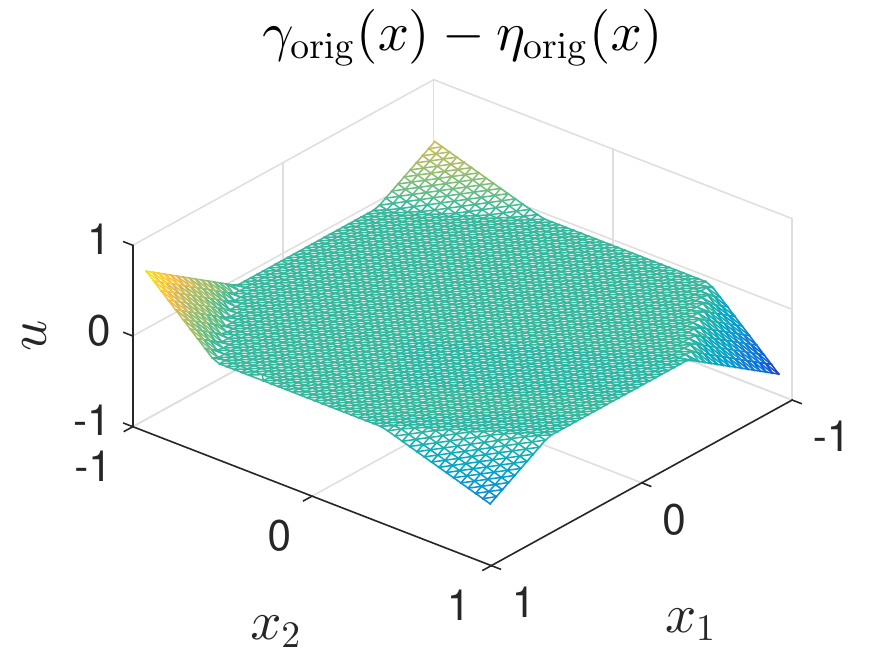}
\caption{PWA explicit MPC law $\mathcal{K}_{\text{orig}}(x)$ (left plot). Decomposition of $\mathcal{K}_{\text{orig}}(x)$ into convex function $\gamma(x)_{\text{orig}} = \mathcal{A}_u^{-1} \circ \gamma(\hat{x})$ and concave function $-\eta_{\text{orig}}(x) = -\mathcal{A}_u^{-1} \circ \eta(\hat{x})$ with $\hat{x} = \mathcal{A}_x \circ x$ (middle plots). The resulting exact representation $\gamma_{\text{orig}}(x) - \eta_{\text{orig}}(x) = \mathcal{A}_u^{-1} \circ (\mathcal{N}(\hat{x};\theta_{\gamma},w,r_{\eta}) - \mathcal{N}(\hat{x};\theta_{\eta},w,r_{\eta}))$ via two deep neural networks is shown on the right plot.}
\label{fig:decomp}
\end{figure*}

We chose the tuning parameters for \eqref{eq:FTOCP} to $P=0$, $R=1$, $Q=2I$ and the horizon to $N=1$. The state constraints are given by $\lvert x_i \rvert \leq 1 \, \text{for} \, i=1,2$ and input constraints by $\lvert u \rvert \leq 1$.
We used the toolbox MPT3 \cite{MPT3} to compute the explicit MPC controller which has 5 regions and is illustrated in the left plot of Fig.~\ref{fig:decomp}.

By applying two invertible affine transformations \eqref{eq:exp_mpc_trans}, the algorithm given in \cite{hempel2013} is used to decompose the explicit MPC controller into the convex function $\gamma_{\text{orig}}(x) = \mathcal{A}_u^{-1} \circ \gamma(\hat{x})$ and the concave function $-\eta_{\text{orig}}(x) = -\mathcal{A}_u^{-1} \circ \eta(\hat{x})$ with $\hat{x} = \mathcal{A}_x \circ x$, depicted in the middle plots of Fig.~\ref{fig:decomp}. 
Both functions consist of $r_\gamma=r_\eta=3$ regions.
According to Theorem~\ref{thm:main_theorem}, two neural networks $\mathcal{N}(\hat{x};\theta_\gamma,3,3)$ and $\mathcal{N}(\hat{x};\theta_\eta,3,3)$ with width $M = n_x + 1 = 3$ and depth $r_\gamma=r_\eta=3$ are used to represent the two convex functions. The parameter values of the networks $\theta_{\gamma}$ and $\theta_{\eta}$ are computed as the minimizers of the mean squared error defined by:
\begin{align}
\theta_\gamma = \underset{\theta_\gamma}{\text{argmin}}\sum_{i = 1}^{n_{\text{tr}}}||\mathcal{N}(\hat{x}_i;\theta_{\gamma}, M, r_\gamma)- \gamma(\hat{x}_i) ||_2^2,
\end{align}
based on $n_{\text{tr}} = 1000$ randomly chosen sampling points for the functions $\gamma(\hat{x})$ (and analogously for $\eta(\hat{x})$).
The learned representation of the neural networks $\gamma_{\text{orig}}(x) - \eta_{\text{orig}}(x) = \mathcal{A}_u^{-1} \circ (\mathcal{N}(\hat{x};\theta_{\gamma},M,r_{\gamma}) - \mathcal{N}(\hat{x};\theta_{\eta},M,r_{\eta}))$ is shown in the right plot of Fig.~\ref{fig:decomp}, which is the same function as the original explicit MPC controller. The training procedure is considered finished when the maximal error $e_{\text{prox}} =\max_{\hat{x}}{|\mathcal{K}(\hat{x}) - (\gamma(\hat{x}) - \eta(\hat{x}))|}$ is less than $0.001$, which we consider to be an exact representation of the transformed explicit MPC law.
The study of the sample complexity of random sampling points $n_{\text{tr}}$ that are necessary to obtain a given error $e_{\text{prox}}$ is an interesting research topic, but it is out of the scope of this paper.
\section{Approximate explicit MPC based on deep learning}\label{sec:approximate}

The previous sections outline two main connections between deep learning and explicit MPC.
The first one is that deep neural networks can exactly represent the explicit MPC law, and not only approximate it arbitrarily well for an increasing number of neurons, as it is known from the universal approximation theorem \cite{barron1993universal}.
The second connection is that, as shown in Lemma~\ref{lemma:expo}, the number of linear regions that deep neural networks can represent grows exponentially with the number of layers. 

While knowing the structure of a network that can exactly represent a given explicit MPC function may be useful in practice, we believe that the use of deep networks to achieve efficient approximations is the most promising idea.
We propose different strategies to deal with the approximation error: a feasibility recovery approach based on control invariant sets and a statistical verification technique to compute safe sets.
Other works have also studied the stability guarantees for approximations of explicit MPC.
In~\cite{kvasnica2011stabilizing} the explicit law is approximated by a polynomial, which lies within a stability tube around the exact MPC law, which guarantees stability.
By using wavelets to approximate the MPC law and barycentric interpolation, \cite{summers2009multiscale} guarantees stability by showing that the cost function of the approximate controller is a Lyapunov function.
If a neural network is used as approximation, giving deterministic stability guarantees is challenging because of the stochastic learning procedure.
In~\cite{hertneck2018learning}, a robust MPC scheme is defined, which accounts for the approximation error of the learning-based approach.
If the assumed error bounds in the robust MPC formulation can be validated a-posteriori by the learned representation, probabilistic statements about stability can be made.
This work focuses on constraint satisfaction, recursive feasibility as well as improved approximation quality while stability guarantees are out of the scope.

\subsection{Training of the deep learning-based approach}

To obtain the proposed approximate explicit MPC law, training data needs to be generated by solving~\eqref{eq:condensed} for $n_{\text{tr}}$ different points.
By randomly choosing $n_{\text{tr,init}}$ initial values and solving~\eqref{eq:condensed} $n_{\text{tr,steps}}$ times in a closed-loop fashion~\eqref{eq:cl_exp} for each initial value, $n_{\text{tr}} = n_{\text{tr,init}} \cdot n_{\text{tr,steps}}$ input samples $x_{\text{tr},i}$ are obtained.
Since from the corresponding optimal control input sequences $\tilde{u}_{\text{tr},i}^*$ only the first step of the sequence is applied to the system in the MPC setting before a new optimal control input sequence is computed, we add $u_0^* = u_{\text{tr},i}$ via the pair $(x_{\text{tr},i},u_{\text{tr},i})$ to the training data set $\mathcal{B}$.
The neural network with a chosen width $M \geq n_x$, a number of layers $L$ is trained with the generated data to find the network parameters $\theta^*$ which minimize the mean squared error over all $n_{\text{tr}}$ training samples:
\begin{align}\label{eq:training}
\underset{\theta}{\text{minimize}}\,\,\,\frac{1}{n_{\text{tr}}}\sum_{i=1}^{n_{\text{tr}}} ||\mathcal{N}(x_{\text{tr},i};\theta, M, L)-u_{\text{tr},i}||^2.
\end{align}
Adam~\cite{kingma2014}, a variant of stochastic gradient descent, is used to solve~\eqref{eq:training} with Keras/Tensorflow \cite{chollet2015keras}, \cite{tensorflow2015-whitepaper}.

The resulting deep network $\mathcal{N}(x;\theta^*,M,L)$ can be used within a feasibility recovery or a verification strategy to control the effect of the approximation error.
These steps are explained in the following subsections.
%
%

\subsection{Feasibility recovery}
Since an exact representation of the explicit MPC law is not achieved, the output of the network is not guaranteed to be a feasible solution of~\eqref{eq:condensed}. In order to guarantee constraint satisfaction as well as recursive feasibility of the problem, the following strategy based on projection is proposed.
The same strategy has been very recently proposed in~\cite{chen2018approximating}.

It is assumed that a convex polytopic control invariant set $\mathcal C_{\text{inv}}$ is available, which is defined as: $\mathcal 
C_{\text{inv}}= \{x\in \mathcal X | \forall x \in \mathcal C_{\text{inv}}, \exists u \in  \mathcal U \text{ s.t. } Ax + Bu \in  C_{\text{inv}}\}$. 
The polytopic control invariant set can be described by a set of linear inequalities as $C_{\text{inv}} = \{ x\in \mathcal X | C_{\text{inv}}x\leq c_{\text{inv}}\}$.

To recover feasibility of the output generated by the neural network, an orthogonal projection  onto a convex set is performed \cite{boyd2004} such that the input constraints are satisfied and the next state lies within the control invariant set:
%
%
\begin{subequations}\label{eq:projection}
\begin{align}
&\underset{\hat{u}}{\text{minimize}} & &\lVert \mathcal{N}(x;\theta,M,L)-\hat{u} \rVert_2^2 \\
&\text{subject to} & &C_{\text{inv}} (A x_{\text{init}} + B \hat u) \leq c_{\text{inv}}, \label{eq:project_states}\\
& & & C_u \hat u \leq c_u,
\end{align}
\end{subequations}
with $x \in \mathcal{C}_{\text{inv}}$ and $\hat{u}^*$ is the optimal and feasible control input to be applied.

Solving~\eqref{eq:projection} directly ensures that the input applied to the system satisfies the input constraints and also 
that the next state satisfies the state constraints, if~\eqref{eq:projection} is feasible.
This in turn means that any consequent state will also belong to $\mathcal C_{\text{inv}}$ because of~\eqref{eq:project_states} and 
therefore problem~\eqref{eq:projection} remains feasible at all times. This also ensures that input and state constraints of the 
closed-loop are satisfied at all times.

\begin{rem}\label{rem:box_constraints}
In the typical case where only box input constraints are present, solving~\eqref{eq:projection} reduces to a saturation operation and no control invariant set is necessary.
In the case of state constraints, a control invariant set should be computed. In the linear case, it is possible to compute such sets even for high dimensional systems \cite{Mirko2017}. The feasibility recovery requires solving the QP \eqref{eq:projection} with $n_u$ variables and $n_{\text{inv}}+n_{\text{cu}}$ constraints, which is often significantly smaller than the original QP defined in~\eqref{eq:condensed} with with $Nn_u$ variables and $N(n_{\text{cx}}+n_{\text{cu}})+n_{\text{cf}}$ constraints. The number of half-spaces $n_{\text{inv}}$ that define a polytopic control invariant set can be reduced if required \cite{Blanco2010} at the cost of conservativeness.
\end{rem}

\begin{rem}
The generation of training points can be simultaneously used for the computation of a control invariant set. For example, if all the vertices of the exact explicit MPC solution are included as training points, taking the convex hull of all of them will generate a control invariant set.
\end{rem}






\subsection{Statistical verification}
The feasibility recovery strategy described in the previous subsection requires the computation of the control invariant set~\cite{Mirko2017}.
But applying this strategy corresponds to solving an additional optimization problem~\eqref{eq:projection} in the control loop when an approximate controller provides an infeasible control input.
Furthermore, it is often the case that the model~\eqref{eq:LTI} used for the design of the controller is just an approximation of the real system which renders the computed control invariant set invalid.
An actual deployment of the controller requires in that case extensive testing on detailed simulators when they are available~\cite{haesaert2015data},~\cite{fan2017d}.

Motivated by this fact, we propose the use of data-driven approaches that enable the a-posteriori statistical verification of the closed-loop performance based on trajectories $\tilde{x} = [x_0,\dots,x_{k_{\text{end}}}]^T$ as explored in~\cite{quendlin_2018}.
The closed-loop dynamics are given by:
\begin{align}\label{eq:cl_dnn}
x_{k+1} = A x_k + B \mathcal{N}(x;\theta,M,L)
\end{align}
for the approximate solution~\eqref{eq:neural_network} and by:
\begin{align}\label{eq:cl_exp}
x_{k+1} = A x_k + B \mathcal{K}(x).
\end{align}
for the explicit MPC controller~\eqref{eq:exp_mpc}.
The data-driven verification can be divided into three steps: data generation, computation of safe sets and validation of safe sets.
The notation of the data sets generated in the verification process is summarized in Table~\ref{tab:description_data_sets} and explained at the corresponding points in this Section.
In this work, we refer to a safe set as the set of initial conditions of a system from which the approximate controller can be applied with a controlled risk of state and input constraint violation.

\begin{table}
\caption{Summary of the accents, subscripts and superscripts used for the notation in the verification for an exemplary data set $\mathcal{A}$ .}
\begin{center}
\begin{tabular}{cc}
Notation & Explanation \\
\midrule
$\tilde{\mathcal{A}}$ & Closed-loop trajectories $\tilde{x}$  \\
\midrule
$\mathcal{A}$         & Initial values $x_0$ of corresponding trajectories in $\tilde{\mathcal{A}}$\\
\midrule
$\tilde{\mathcal{A}}^+$ & Closed-loop trajectories for which no violations occured \\
\midrule
$\tilde{\mathcal{A}}^-$ & Closed-loop rajectories for which violations occured \\
\midrule
$\tilde{\mathcal{A}}_{\text{exp}}$ & Closed-loop trajectories generated via exact MPC~\eqref{eq:cl_exp} \\
\midrule
$\tilde{\mathcal{A}}_{\text{dnn}}$ & Closed-loop trajectories generated via NN controller~\eqref{eq:cl_dnn} \\
\bottomrule
\end{tabular}
\label{tab:description_data_sets}
\end{center}
\end{table}


\subsubsection{Data generation}

For the verification procedure, different data sets are necessary, which are all independent from the data used for training the neural networks in~\eqref{eq:training}.
The data sets needed for verification are $\tilde{\mathcal{D}} \in \{ \tilde{\mathcal{G}}_{\text{dnn}}, \tilde{\mathcal{T}}_{\text{dnn}}, \tilde{\mathcal{T}}_{\text{exp}}, \tilde{\mathcal{V}}_{\text{dnn}} \}$, where the \emph{\textasciitilde} denotes that they contain closed-loop trajectories $\tilde{x}$.
The set $\tilde{\mathcal{G}}_{\text{dnn}}$ is used to compute the safe set with two different methods.
The sets $\tilde{\mathcal{T}}_{\text{dnn}}$ and $\tilde{\mathcal{T}}_{\text{exp}}$ are used to compare the size of the safe sets for the approximate and the exact controller.
The initial values of the corresponding trajectories in $\tilde{\mathcal{T}}_{\text{dnn}}$ and $\tilde{\mathcal{T}}_{\text{exp}}$ are identical.
The validation set $\tilde{\mathcal{V}}_{\text{dnn}}$ is used to compute the probability with which a trajectory starting from an initial value within a computed safe set does not violate the constraints.
Sets denoted with a \emph{dnn} in the subscript contain trajectories obtained via closed-loop simulation with the approximate controller~\eqref{eq:cl_dnn} whereas the subscript \emph{exp} indicates trajectories obtained with the exact MPC~\eqref{eq:cl_exp}.
The data sets are evaluated by requirements formulated using metric temporal logic (MTL) \cite{baier2008principles}:
\begin{align} 
\rho(\tilde{x}) = \square_{\interval{0}{k_{\text{end}}}}(C_x x_k \leq c_x \land C_u u_k \leq c_u), \label{eq:log_requirement}
\end{align}
where $\land$ is the operator for \emph{logical and} and $C_x$ and $c_x$ and $C_u$ and $c_u$ describe the state and input constraints as defined in~\eqref{eq:FTOCP}.
This requirement translates to $C_x x_k$ and $C_u u_k$ have to be \emph{always} smaller than or equal to $c_x$ and $c_u$ between time steps $0$ and $k_{\text{end}}$ for the requirements to be satisfied.
If this is true, we obtain $\rho(\tilde{x}) = +1$ and $\rho(\tilde{x}) = -1$ otherwise.
By replacing $\square$ with $\lozenge$ in \eqref{eq:log_requirement} the meaning would change from \emph{always} to \emph{eventually}.
More complex requirements can be included using MTL~\cite{baier2008principles},~\cite{quendlin_2018}.

The requirement~\eqref{eq:log_requirement} is used to evaluate the previously mentioned sets of trajectories $\tilde{\mathcal{D}} \in \{ \tilde{\mathcal{G}}_{\text{dnn}}, \tilde{\mathcal{T}}_{\text{dnn}}, \tilde{\mathcal{T}}_{\text{exp}}, \tilde{\mathcal{V}}_{\text{dnn}} \}$.
Thus, each data set is split into one containing all valid trajectories (marked with \emph{$+$} in the exponent) and into one containing all invalid trajectories (marked with \emph{$-$} in the exponent), for instance $\tilde{\mathcal{D}}^+ = \{\tilde{x} \, | \, \tilde{x} \in \tilde{\mathcal{D}} \land \rho(\tilde{x}) = +1 \}$.
The set $\mathcal{D} \in \{ \mathcal{G}_{\text{dnn}}, \mathcal{T}_{\text{dnn}}, \mathcal{T}_{\text{exp}}, \mathcal{V}_{\text{dnn}} \}$ contains the initial value of every trajectory in the corresponding sets $\tilde{\mathcal{D}}$.


\subsubsection{Safe sets}

Two approaches are applied to obtain an explicit description of the safe set from closed-loop data containing factual negatives and positives.
After deriving a probabilistic safe set $\mathcal{S}$, we define all $x \in \mathcal{S}$ as probabilistic positives, and all $x \not\in \mathcal{S}$ as probabilistic negatives.
The first approach uses an $m$-dimensional ellipsoidal given by $x^T E x = 1$ with $E \in \mathbb{R}^{m \times m}$ and $E = E^T$.
With the following convex optimization problem, it is possible to find the ellipsoidal safe set that inscribes the maximum-sized hypercube such that no initial values of invalid trajectories are contained in said ellipsoidal safe set:
\begin{subequations}\label{eq:ellipsoid}
	\begin{align}
		&\underset{E}{\text{minimize}} && \Tr(E) \span \span \span \\
		&\text{subject to} && E \succeq 0, \\
		& && x^T E x \geq (1+\epsilon), \,\, \forall x \in \mathcal{G}_{\text{dnn}}^-
	\end{align}
\end{subequations}
where $\epsilon \geq 0$ is a tuning parameter to increase robustness, which is desirable due to the finite number of points in $\mathcal{G}_{\text{dnn}}^-$ used for the computation of the ellipsoidal.
The resulting ellipsoidal set is then described by
\begin{align}\label{eq:ellipsoidal_safe_set}
\mathcal{S}_{\text{ell}} = \left\{ \, x \, \middle\vert \, x^T E x \leq 1\right\}.
\end{align}
The value of $\epsilon$ was in this work tuned via trial and error, such that $x^T E x > 1$ for all $x \in V^-_{\text{dnn}}$, which is equivalent to the absence of false positives in the validation set.

The second approach relies on support vector machines (SVMs) for classification~\cite{cortes1995support} to derive a less restrictive safe set.
The SVM learning problem is given by:
\begin{subequations}\label{eq:SVM}
	\begin{align}
		&\underset{w,v,\zeta}{\text{minimize}} && \frac{1}{2}w^Tw + C \sum_{m=1}^{|\mathcal{G}_{\text{dnn}}|}{\zeta_m} \span \span \span \\
		&\text{subject to} && y_m(w^T \phi(x_m)+v) \geq 1- \zeta_m, \\
		& && \zeta_m \geq 0, \\
		& && x_m \in \mathcal{G}_{\text{dnn}},\\		
		& && \forall m = 1,\dots,|\mathcal{G}_{\text{dnn}}|, 
	\end{align}
\end{subequations}
where $w$ is the weight vector, $v$ is the bias and $\phi(\cdot)$ defines a kernel $\kappa(x_i,x_j) = \phi(x_i)^T\phi(x_j)$  and
$\zeta_m$ are slack variables to relax the problem. The tuning parameter $C>0$ is a penalty term to weigh the importance of misclassification errors and $y_m$ are the decision functions.
The resulting relaxed safe set will only be dependent on the $m_{\text{cr}}$ decision functions $y_i$ where $\zeta_i = 0$.
The corresponding $x_i$ are called support vectors and define the set description:
\begin{align}\label{eq:SVM_safe_set}
\mathcal{S}_{\text{SVM}} = \{x \, | \, y_i(w^T \kappa(x,x_i)+v) \geq 1, \, \forall i = 1,\dots,m_{\text{cr}}\}.
\end{align}

\begin{figure}[t]
\begin{center}
\includegraphics[width=0.99\columnwidth]{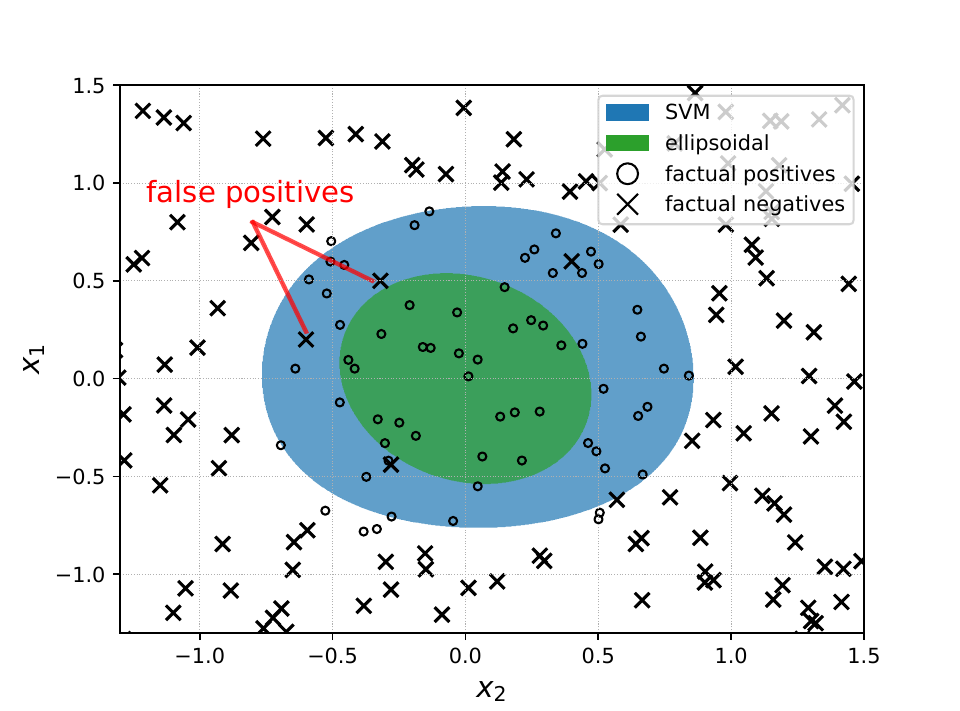} 
\caption{Exemplary visualization of the safe sets obtained from data containing factual positives (circles) and factual negatives (marks) for a system with two states $x_1$ and $x_2$. The ellipsoidal robust safe set contains only factual positives, but the less conservative safe set computed by a SVM contains factual negatives which are considered as false positives.}
\label{fig:safe_sets}
\end{center}
\end{figure}

The two types of safe sets are illustrated exemplarily in Fig.~\ref{fig:safe_sets}.
The circles depict initial values of trajectories $x = \tilde{x}(0)$ without constraint violations ($\rho(\tilde{x})=+1$, factual positives) whereas crosses represent those of invalid trajectories ($\rho(\tilde{x})=-1$, factual negatives). 
The ellipsoidal set (green) is more conservative, but it does not contain any initial values that lead to constraint violations.
The larger set (blue) computed by SVMs allows a controlled trade-off between the size of the safe set and the proneness towards misclassification.

Two different types of classification errors can be distinguished taking the view of the probabilistic safe sets.
False negatives are initial values from which trajectories are falsely classified by a safe set as leading to constraint violations whereas false positives are initial values from which trajectories are considered to be safe but in fact lead to constraint violations.
False positives are the worst case since they might lead to actual constraint violations of the closed-loop system.

In order to assess the conservativeness of a computed safe set $\mathcal{S} \in \{ \mathcal{S}_{\text{ell}}, \mathcal{S}_{\text{SVM}} \}$, its size is compared to the size of the safe set that is obtained with the exact controller.
The initial values in the test data set $\mathcal{T}_{\text{dnn}}$ are evaluated via~\eqref{eq:ellipsoidal_safe_set} or \eqref{eq:SVM_safe_set} to obtain the set $\mathcal{S}^+_{\mathcal{T}_{\text{dnn}}} = \{ x \, | \, x \in \mathcal{S} \land x \in \mathcal{T}_{\text{dnn}} \}$, containing all initial values which are part of the safe set.
We define the volume of the computed safe set in relation to the exact one as:
\begin{align}
m_{\text{vol}} = \frac{|\mathcal{S}^+_{\mathcal{T}_{\text{dnn}}} \cap \mathcal{T}_{\text{exp}}^+|}{|\mathcal{T}_{\text{exp}}^+|}, \label{eq:m_vol}
\end{align}
where the numerator gives the number of initial values in the test set which lead to feasible trajectories for the approximate and exact solution and the denominator gives the number of all initial values which lead to feasible trajectories with the exact solution.


\subsubsection{Validation}

Since the computation of the safe sets is based on a finite amount of data points, we analyze the approximation quality of the computed safe sets by computing the amount of false positives as well as by applying methods from statistical learning theory (STL)~\cite{luxburg2011statistical}.
The proportion of false positives for classifying the initial values in $\mathcal{V}_{\text{dnn}}$ is given by:
\begin{align}
m_{\text{fp}} = \frac{|\mathcal{S}^+_{\mathcal{V}_{\text{dnn}}} \cap \mathcal{V}_{\text{dnn}}^-|}{|\mathcal{S}^+_{\mathcal{V}_{\text{dnn}}}|}, \label{eq:m_fp}
\end{align}
where $\mathcal{S}^+_{\mathcal{V}_{\text{dnn}}} = \{ x \, | \, x \in \mathcal{S} \land x \in \mathcal{V}_{\text{dnn}} \}$.

To make probabilistic statements about the safe set with a certain confidence, STL is used.
While~\cite{hertneck2018learning} recently used STL to define an upper bound on the approximation error of an approximate controller to make statements about the closed-loop behavior, in this work STL is used to directly validate the closed-loop performance of the approximate solution. This means that the validation strategy is valid regardless of the (potentially incorrect) model that is used to compute the MPC solution provided that a simulator of the real system exists.
An indicator function $I(x)$ is introduced to assign to all initial values a risk via the temporal logic requirement defined in~\eqref{eq:log_requirement}:
\begin{align}\label{eq:indicator}
I(x) = 
\begin{cases}
1 & \text{if} \, \, \rho(\tilde{x}) = +1,\\
0 & \text{if} \, \, \rho(\tilde{x}) = -1,
\end{cases}
\end{align}
with $x = \tilde{x}(0)$.
The expected value of \eqref{eq:indicator} for all $x \in \mathcal{S}_{\mathcal{V}_{\text{dnn}}}^+$ describes the empirical risk $r_{\text{emp}}$ that a trajectory starting from initial value $x \in \mathcal{S}$ leads to constraint satisfaction over the whole time period $[0, k_{\text{end}}]$:
\begin{align}\label{eq:empirical_risk}
r_{\text{emp}} = \frac{1}{n_{\text{s}}}\sum_{i=1}^{n_{\text{s}}}{I(x_i)}
\end{align}
where $n_{\text{s}} = |\mathcal{S}_{\mathcal{V}_{\text{dnn}}}^+|$.
The empirical risk $r_{\text{emp}}$ is based on a data set and is the best approximation of the true risk $r_{\text{true}}$, which is based on all $x \in \mathcal{S}$.
Hoeffdings inequality \cite{hoeffding1994probability} provides an upper bound on the probability that the empirical risk $r_{\text{emp}}$ deviates more than $\delta$ from the true risk $r_{\text{true}}$:
\begin{align}\label{eq:hoeffding}
P \left(  \left| r_{\text{emp}} - r_{\text{true}} \right| \geq \delta\right) \leq 2 \exp(-2n_{\text{s}}\delta^2).
\end{align}
Thus, for all $x \in \mathcal{S}$
\begin{align}\label{eq:safety_and_confidence}
P(I(x)=1) = r_{\text{true}} \geq r_{\text{emp}} - \delta,
\end{align}
with confidence level $h_{\delta} = 1 - 2 \exp(-2n\delta^2)$.
The meaning of~\eqref{eq:safety_and_confidence} is that a trajectory starting from an initial value within the safe set will not violate the constraints with a probability greater than or equal to $r_{\text{emp}} - \delta$.

\subsection{Alternative approximation methods}



The proposed deep learning-based approximate explicit MPC approach is compared to other approximation approaches.

The first alternative approximates the explicit controller using multi-variate polynomials of the form $\mathcal{P}:\mathbb{R}^{n_x} \rightarrow \mathbb{R}^{n_u}$ with degree $p$:

\begin{align}\label{eq:approx_poly}
\mathcal{P}(x;\alpha,p) =
\begin{bmatrix}
\overset{p}{\underset{i_{1}=0}{\sum}} \dots \overset{p}{\underset{i_{n_x}=0}{\sum}} a_{1,m} \prod_{j=1}^{n_x} x_j^{i_{j}} \\
\vdots \\
\overset{p}{\underset{i_{1}=0}{\sum}} \dots \overset{p}{\underset{i_{{n_x}}=0}{\sum}} a_{n_u,m} \prod_{j=1}^{n_x} x_j^{i_{j}} \\
\end{bmatrix}
\end{align}
where the index $m = \sum_{j=1}^{n_x} i_j$ and $\alpha_i = \{a_{i,1},\dots,a_{i,(p+1)^{n_x}} \}$ for $i=1,\dots,n_u$ contains all coefficients. The coefficients of the polynomials are computed by solving
\begin{align}\label{eq:training_poly}
\underset{\alpha}{\text{minimize}}\,\,\,\frac{1}{n_{\text{tr}}}\sum_{i=1}^{n_{\text{tr}}} ||\mathcal{P}(x_{\text{tr},i};\alpha,p)-u_{\text{tr},i}||^2
\end{align}
where $u_{\text{tr},i}$ is the exact optimal control input obtained solving~\eqref{eq:condensed} for each training point $x_{\text{tr},i}$.
The memory footprint of a multi-variate polynomial is given by
\begin{equation}
\Gamma_{\mathcal{P}} = \alpha_{\text{bit}} n_u (p+1)^{n_x}.
\end{equation}

The second method is similar to the approach in \cite{holaza2013}.
We use the partition of an explicit MPC description with a shorter horizon $N \leq N_{\text{max}}$ (and therefore less regions) and adapt the parameters $\lambda = \{\lambda_1,\dots,\lambda_{n_{\text{r}}}\}$ where $\lambda_i = \{K_i,g_i\}$ by solving the following optimization problem:
\begin{align}\label{eq:training_optimized}
\underset{\lambda}{\text{minimize}}\,\,\,\frac{1}{n_{\text{tr}}}\sum_{i=1}^{n_{\text{tr}}} ||\mathcal{L}_N(x_{\text{tr},i};\lambda)-u_{\text{tr},i}||^2.
\end{align}
We denote  the optimized descriptions $\mathcal{L}:\mathbb{R}^{n_x} \rightarrow \mathbb{R}^{n_u}$ as
\begin{align}\label{eq:approx_smaller_horizon}
\mathcal{L}_N(x;\lambda) = 
	\begin{cases}
		K_1 x+g_1 & \text{if} \quad x \in \mathcal{R}_1, \\
		& \vdots \\
		K_r x+g_r & \text{if} \quad x \in \mathcal{R}_r.
	\end{cases}
\end{align}
The memory footprint of the optimized explicit MPC can be estimated as done for the standard explicit MPC~\eqref{eq:memory_empc}.



%

The explicit MPC description and the approximation methods are summarized in Table~\ref{tab:parameters}. We introduce a new abbreviation for explicit MPC laws $\mathcal K_N(x)$ where $N$ stands for the horizon of the primary problem \eqref{eq:FTOCP} they are derived from.

\begin{table}
\caption{Summary of algorithms used, including exact explicit solution $\mathcal K_N$ and the approximations methods.}
\begin{center}
\begin{tabular}{lcc}
Method & Param. & Explanation \\
\midrule
$\mathcal{K}_N(x)$ & $N$ & prediction horizon  \\
\hline 
\multirow{3}{*}{$\mathcal{N}(x;\theta,M,L)$} & $\theta$ & aff. trans. $\{W_l,b_l\}$ $\forall$ layers \\ 
& $M$ & neurons per hidden layer \\ 
& $L$ &number of hidden layers\\ 
\hline 
\multirow{2}{*}{$\mathcal{P}(x;\alpha,p)$} & $p$ & degree of the polynomial \\ 
& $\alpha$ & coefficients $a_i$ for all terms\\ 
\hline
\multirow{2}{*}{$\mathcal{L}_N(x;\lambda)$} & $N$ & prediction horizon \\
& $\lambda$ & aff. trans. $\{K_i,g_i\}$ $\forall$ regions \\
\bottomrule
\end{tabular}
\label{tab:parameters}
\end{center}
\end{table}
\section{Simulation results}\label{sec:results}
The potential of the proposed approach is illustrated with a simulation example modified from \cite{wang2010} and the classic example of the inverted pendulum on a cart.
The goal is in both cases to steer the system to the origin.

The approximate methods via polynomials $\mathcal{P}(x;\alpha,p)$~\eqref{eq:approx_poly}, optimized explicit MPC with reduced horizon $\mathcal{L}(x;\lambda)$~\eqref{eq:approx_smaller_horizon} and neural networks $\mathcal{N}(x;\theta,M,L)$~\eqref{eq:neural_network} are compared with respect to their performance and their memory footprint~\eqref{eq:memory_empc} to the exact explicit MPC solution $\mathcal{K}(x)$~\eqref{eq:exp_mpc}.
The exact explicit MPC controller is considered as the benchmark for the chosen performance index average settling time (AST). 
The AST is defined as the time necessary to steer all states to the origin.
A state is considered to be at the origin when $\lvert x_i \rvert \leq \SI{1e-2}{}$.
The relative AST (rAST) is the performance measure with respect to the exact solution with the longest horizon $N_{\text{max}}$.

In the following, the dependency of the controllers on $x$ and on the parameters are dropped for the sake of brevity.
Additionally, neural networks $\mathcal{N}(x;\theta,M,L)$ will be referred to by $\mathcal{N}_{M,L}$ and polynomials $\mathcal{P}(x;\alpha,p)$ by $\mathcal{P}_p$.

\begin{rem}
We investigate in this section the performance of shallow ($L=1$) and deep ($L\geq2$) neural networks. For very deep networks ($L \gg 10$), the vanishing gradient problem can occur in the training phase which jeopardizes the approximation accuracy. To counteract the effect, measures like highway layers~\cite{srivastava2015highway} can be taken. In this work, applying countermeasures is not necessary, because the used ReLU networks are less prone to vanishing gradients and the deepest network considered does not exceed $L=10$ layers.
\end{rem}



\subsection{Case-studies}
Two examples to investigate the proposed approach are introduced.
The control tasks are solved many times from different initial conditions.
The trajectories generated with the corresponding explicit MPC solutions $\mathcal{K}_7$ and $\mathcal{K}_{10}$ were used to train the different approximation approaches \eqref{eq:training}, \eqref{eq:training_poly} and \eqref{eq:training_optimized}.
Since both case-studies include box input constraints, a simple saturation was used to guarantee satisfaction of the input constraints for the approximate controller as proposed in Remark~\ref{rem:box_constraints}.
It is assumed that all states of the systems can be measured.

\subsubsection{Oscillating Masses (OM)}
\begin{figure}[t]
\begin{center}
\includegraphics[width=0.7\columnwidth]{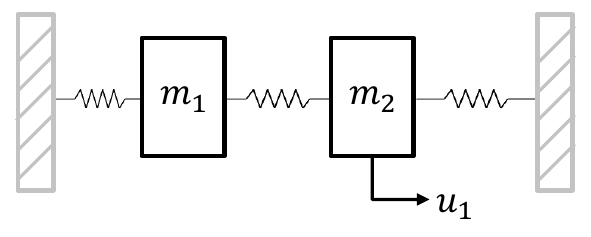} 
\caption{Chain of masses connected via springs.}
\label{fig:oscillating_masses}
\end{center}
\end{figure}

The first example represents two horizontally oscillating masses interconnected via a spring where each one is connected via a spring to a wall, as shown in Fig.~\ref{fig:oscillating_masses}.
Both masses can only move horizontally and have a weight of \SI{1}{\kilogram} and each spring has a constant of \SI{1}{\newton\per\metre}. The states of each mass are its position, limited to $\lvert s \rvert \leq \SI{4}{\metre}$, and its speed $v$, limited to $\lvert v \rvert \leq \SI{10}{\metre\per\second}$.
A force limited by $\lvert u \rvert \leq \SI{.5}{\newton}$ can be applied to the right mass.

The state vector is given by $x = [s_1, v_1, s_2, v_2]^T$ and the system matrices are discretized with first-order hold and a sampling time of \SI{0.5}{\second} resulting in:
 
\begin{equation*}
\begin{aligned}
&A=
\begin{bmatrix}
0.763   & 0.460 & 0.115 & 0.020 \\ 
-0.899  & 0.763 & 0.420 & 0.115 \\
0.115  & 0.020 & 0.763 & 0.460 \\
0.420  & 0.115 & -0.899 & 0.763 \\
\end{bmatrix},
&B=
\begin{bmatrix}
0.014 \\
0.063 \\
0.221 \\
0.367 \\
\end{bmatrix}.
\end{aligned}
\end{equation*}



The benchmark horizon was $N_{\text{max}}=7$ corresponding to 2317 regions.
The exact explicit controller $\mathcal{K}_7$ was used to generate 25952 training samples .


\subsubsection{Inverted pendulum on cart (IP)}

\begin{figure}[t]
\begin{center}
\includegraphics[width=0.9\columnwidth]{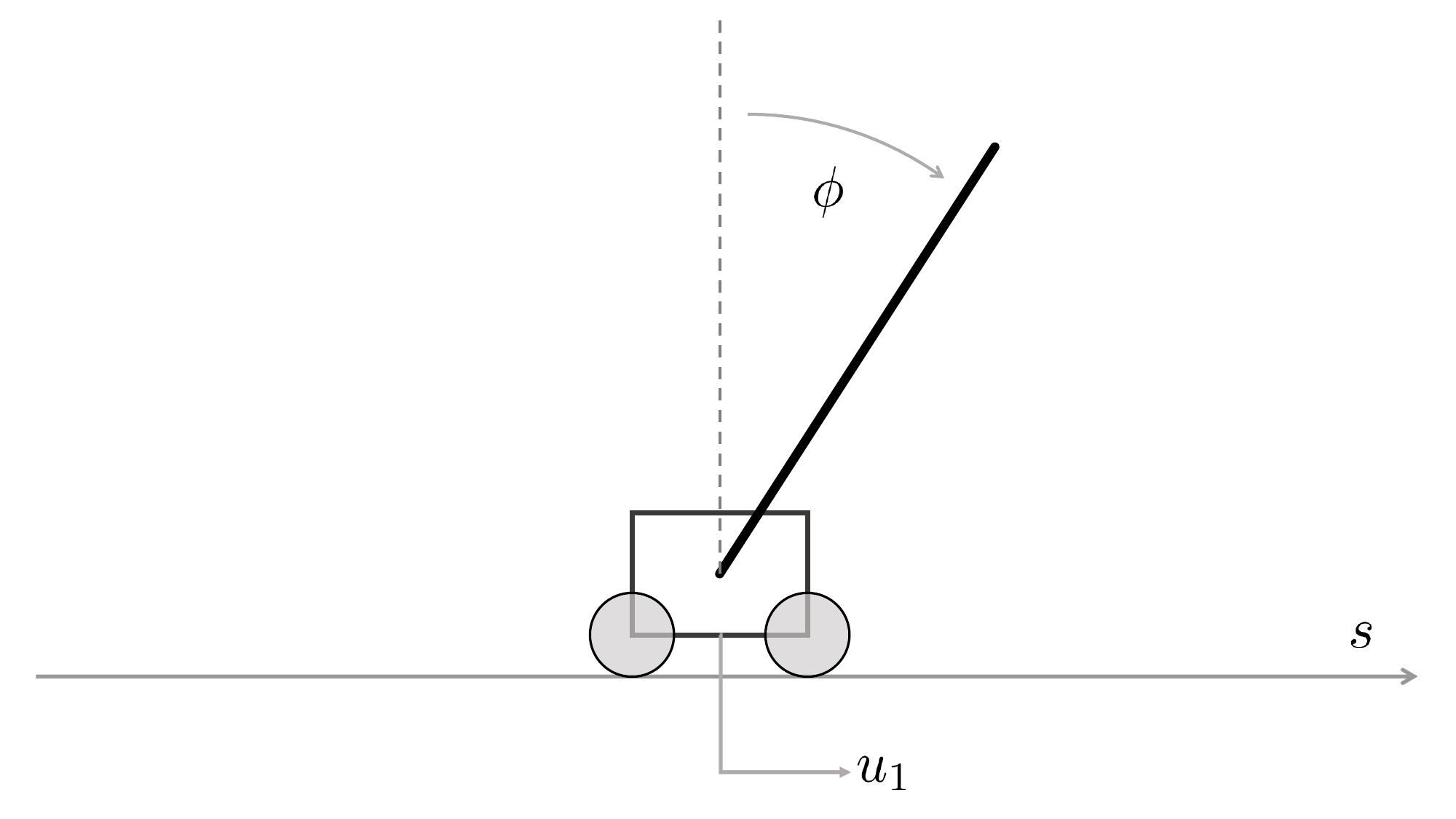}
\caption{Inverted pendulum on a cart.}
\label{fig:inverted_pendulum_on_cart}
\end{center}
\end{figure}

The second example is the inverted pendulum on a cart, illustrated in Fig.~\ref{fig:inverted_pendulum_on_cart}.
The goal is to keep the pole erected and the cart in the central position.
The states are the angle of the pole $\Phi$, its angular speed $\dot{\Phi}$, the position of the cart $s$ and the speed of the cart $\dot{s}$.
The states $x = [\Phi,s,\dot{\Phi},\dot{s}]^T$ are constrained to $|x|^T \leq [1, 1.5, 0.35, 1.0]^T$.
The force $|u| \leq \SI{1}{\newton}$ is directly applied to the cart.

Euler-discretization with a sampling time of \SI{0.1}{\second} was used to obtain the discrete system dynamics given by:
\begin{equation*}
\begin{aligned}
&A=
\begin{bmatrix}
1   & 0.1 & 0 & 0 \\ 
0 & 0.9818 & 0.2673 & 0 \\
0  & 0  & 1 & 0.1 \\
0 & -0.0455 & 3.1182 & 1 \\
\end{bmatrix},
&B=
\begin{bmatrix}
0 \\
0.1818 \\
0 \\
0.4546 \\
\end{bmatrix}.
\end{aligned}
\end{equation*}

For this example the explicit benchmark solution was computed with horizon $N_{\text{max}}=10$ resulting in a PWA function consisting of 1638 polyhedral regions.
88341 samples were generated to train the approximated controllers.

\subsection{Performance}

We investigated both examples OM and IP by simulating closed-loop trajectories starting from randomly chosen initial values within the feasible state space.
For each initial value, the exact MPC controller $\mathcal{K}_{(\cdot)}$ and the approximate methods were applied.
The controller $\mathcal{K}_{(\cdot)}$ provided the benchmark performance.
For both case-studies, the evaluation led to similar results, as it can be seen in Table~\ref{tab:AST_and_memory}.
The proposed deep neural networks $\mathcal{N}_{6,6}$ and $\mathcal{N}_{10,6}$ only use $\SI{0.23}{\percent}$ and $\SI{1.07}{\percent}$ of the memory of the optimal solutions $\mathcal{K}_7$ and $\mathcal{K}_{10}$ while reaching an average AST that is only $\SI{1.5}{\percent}$ and $\SI{3.8}{\percent}$ longer than the exact solution. The deep neural network clearly achieves the best trade-off between performance and memory requirements. It is interesting to see that a shallow network $\mathcal{N}_{43,1}$ and $\mathcal{N}_{120,1}$ with a slightly larger memory footprint than the deep network, achieves considerably worse performance.
The results show that a naive polynomial approximation of the explicit MPC does not lead to good results as the performance that can be achieved with no more than 2 $\SI{}{\kilo\byte}$ is significantly worse than the other approximation methods.
Even if the optimized explicit MPC with the finest partition $\mathcal{L}_6$ and $\mathcal{L}_7$ is compared to the benchmark, the proposed deep neural network performs slightly better for OM and clearly better for IP while having a much smaller memory footprint.


\begin{figure}[t]
\begin{center}
\includegraphics[width=0.99\columnwidth]{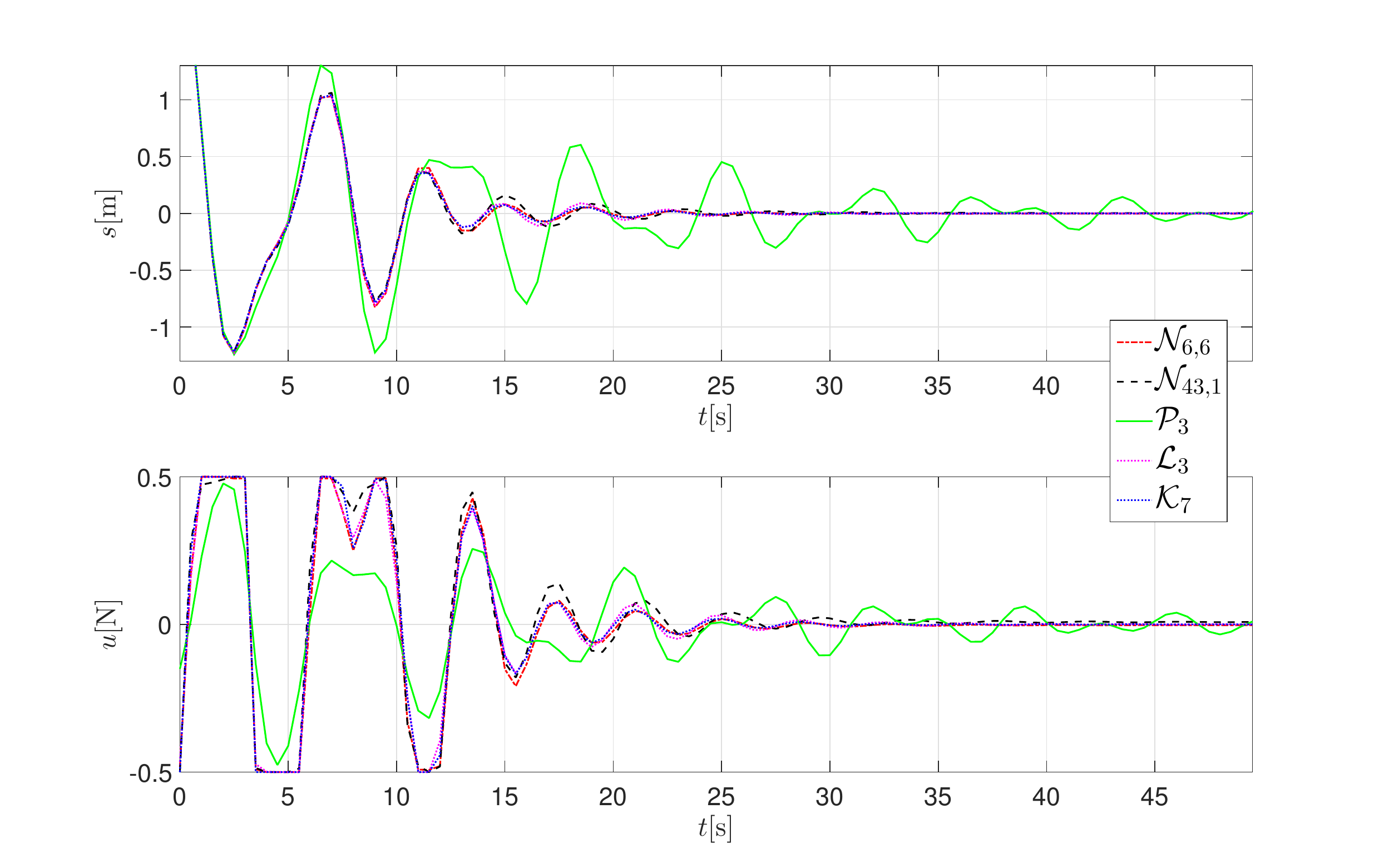}
\caption{Position of the first mass (top plot) and control inputs (bottom plot) for different control strategies for one exemplary closed-loop simulation of the oscillating masses.}
\label{fig:results}
\end{center}
\end{figure}

\begin{figure}[t]
\begin{center}
\includegraphics[width=0.99\columnwidth]{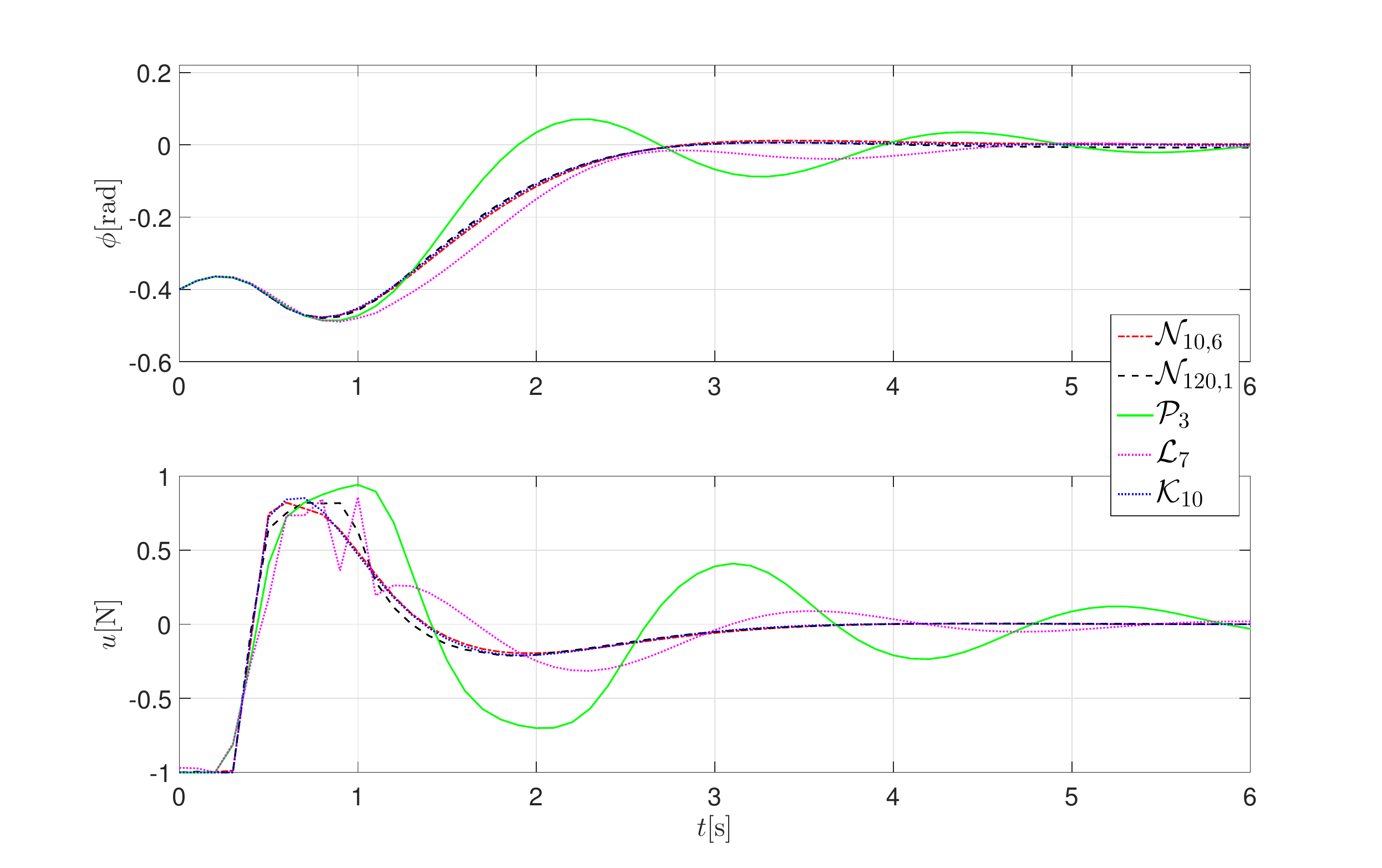}
\caption{Position of the pendulum (top plot) and control input (bottom plot) for different control strategies for one exemplary closed-loop simulation of the inverted pendulum on a cart.}
\label{fig:performance_inverted_pendulum_on_cart}
\end{center}
\end{figure}

Fig.~\ref{fig:results} and Fig.~\ref{fig:performance_inverted_pendulum_on_cart} show an example of the closed-loop trajectories obtained for each type of controller for the two examples. It can be clearly seen that the polynomial approximation (degree 3) cannot properly approximate the explicit controller. The best results, which are almost identical to the exact explicit controller $\mathcal K_7$ and $\mathcal K_{10}$, are obtained by the deep neural networks $\mathcal{N}_{6,6}$ and $\mathcal{N}_{10,6}$.

\begin{table}
	\centering
	\caption{Comparison of the relative average settling time (rAST) for 10000 simulation runs and memory footprint $\Gamma$ for different controllers for the oscillating masses (OM) and inverted pendulum on cart (IP) example.}
	\begin{tabular}{ccccccc}
	&&&&&&\\
OM	& $\mathcal{K}_7$ & $\mathcal{L}_6$ & $\mathcal{L}_3$ & $\mathcal{P}_3$ & $\mathcal{N}_{6,6}$ & $\mathcal{N}_{43,1}$ \\
	\midrule
	rAST  [-]& 1  & 1.020 & 1.113 & 1.407 & 1.015 & 1.125\\
	$\Gamma \,[\SI{}{\kilo\byte}]$ & 691.9   &  431.8 & 38.3  & 2.00  & 1.93  & 2.02 \\
	\bottomrule
	&&&&&&\\
	&&&&&&\\
IP	& $\mathcal{K}_{10}$ & $\mathcal{L}_7$ & $\mathcal{L}_6$ & $\mathcal{P}_3$ & $\mathcal{N}_{10,6}$ & $\mathcal{N}_{120,1}$ \\
	\midrule
	rAST  [-]& 1  & 1.897 & 2.273 & 2.276 & 1.038 & 1.060\\
	$\Gamma \,[\SI{}{\kilo\byte}]$ & 444.7   &  191.5 & 137.1  & 2.00  & 4.77  & 5.63 \\
	\bottomrule
	\end{tabular}
	\label{tab:AST_and_memory}
\end{table}


\subsection{Statistical verification}
To derive the set in which the application of the proposed controller is safe, explicit descriptions of safe sets are computed via~\eqref{eq:ellipsoid}~and~\eqref{eq:SVM}.
The data sets containing the initial values of the trajectories and their cardinality are
$|\mathcal{G}_{\text{dnn}}| = 20000$, 
$|\mathcal{T}_{\text{exp}}| = |\mathcal{T}_{\text{dnn}}| = 10000$ and 
$|\mathcal{V}_{\text{dnn}}| = 40000$.
In case of an exact approximation $\mathcal{K}_7$ and $\mathcal{K}_{10}$ and $\mathcal{N}_{6,6}$ and $\mathcal{N}_{10,6}$ would be equivalent.
This would mean that all initial values $x \in \mathcal{D}^+$ are part of the safe set.
Since $\mathcal{N}_{6,6}$ and $\mathcal{N}_{10,6}$ are only approximations, the test sets $\mathcal{T}_{\text{dnn}}^+$ and $\mathcal{T}_{\text{exp}}^+$ are directly compared to derive a first naive measure of the approximation quality.
The ratio of the cardinality of the two sets:
\begin{align}
m_{\text{dir}} = \frac{|\mathcal{T}_{\text{dnn}}^+|}{|\mathcal{T}_{\text{exp}}^+|}, \label{eq:m_dir}
\end{align}
is \SI{97.7}{\percent} for OM and \SI{95.1}{\percent} for IP as pointed out in the rows of Table~\ref{tab:safe_set_and_approximation} denoted by \emph{direct}.

\begin{table}
	\centering
	\caption{Comparison of approximate neural network controllers and exact explicit MPC and their safe sets for the examples oscillating masses (OM) and inverted pendulum on cart (IP). }
	\begin{tabular}{ccccc}
	OM & vol. [\SI{}{\percent}] & false pos. [\SI{}{\percent}] & safety [\SI{}{\percent}] & confidence [\SI{}{\percent}] \\
	\midrule
	direct        & 97.7 &  -    &  -     &     -    \\
	ellipsoidal   & 70.5 &  0    &  98.5  &  >99.9  \\
	SVM           & 98.3 &  1.6  &  96.9  &  >99.9  \\
	\bottomrule
	&&&& \\
	&&&& \\
	IP & vol. [\SI{}{\percent}] & false pos. [\SI{}{\percent}] & safety [\SI{}{\percent}] & confidence [\SI{}{\percent}] \\
	\midrule
	direct        & 95.1 &  -   &  -   &     -    \\
	ellipsoidal   & 53.3 &  0   & 97.0 &   >99.9  \\
	SVM           & 83.1 & 1.8  & 95.1 &   >99.9  \\
	\bottomrule
	\end{tabular}
	\label{tab:safe_set_and_approximation}
\end{table}

The optimization problem \eqref{eq:ellipsoid} is solved to obtain explicit formulations of an ellipsoidal safe set $\mathcal{S}_{\text{ell}}$.
The estimated volume~\eqref{eq:m_vol} for OM is \SI{70.5}{\percent} and \SI{53.3}{\percent} for IP. 
This restricts the usage of the approximated controller to a significantly smaller volume, but the absence of false positives on the validation set ($m_{\text{fp}} = 0$)~\eqref{eq:m_fp} indicates a certain robustness, only limited in this case by the use of a finite amount of data points.
By applying~\eqref{eq:safety_and_confidence} for both examples it can be said with confidence $>\SI{99.9}{\percent}$ that a trajectory starting from $x \in \mathcal{S}_{\text{ell}}$ will not violate the constraints with a probability $\geq\SI{97.0}{\percent}$.
The results considering $\mathcal{S}_{\text{ell}}$ are given in the rows of Table~\ref{tab:safe_set_and_approximation} named \emph{ellipsoidal}.

To overcome the conservativeness of $\mathcal{S}_{\text{ell}}$ with respect to the covered volume $m_{\text{vol}}$, a relaxed safe set $\mathcal{S}_{\text{SVM}}$ is computed via~\eqref{eq:SVM}.
Since the data is not linearly separable the radial basis function (RBF) is chosen as the kernel:
\begin{align}
\kappa(x, x_i) = \text{exp}\left( -\nu \mid\mid x-x_i \mid\mid^2\right),\label{eq:rbf}
\end{align}
where $\nu$ is a tuning parameter defining the width of the kernel function.
The volume of the SVM safe sets are $\SI{98.3}{\percent}$ for OM and $\SI{83.1}{\percent}$ for IP while classifying less than $\SI{2}{\percent}$ of the validation set as false positives.
This proportion of false positives leads to a slightly reduced safety of $\SI{96.9}{\percent}$ for OM and $\SI{95.1}{\percent}$ for IP with confidence $>\SI{99.9}{\percent}$
The results are summarized in the rows of Table~\ref{tab:safe_set_and_approximation} labelled \emph{SVM}.

This shows that the proposed controller can be applied within the safe sets with a small risk of constraint violation with high confidence.
Especially if the application can be allowed in $\mathcal{S}_{\text{SVM}}$, the volume of the safe set is similar to the explicit MPC solution while providing a comparable performance.
For safety critical applications, fallback strategies should be available to be used when the approximate controller cannot achieve the desired performance.

\subsection{Embedded implementation}

The embedded implementation of the proposed approximate neural network controllers is straightforward since evaluating neural networks consists only of multiplications, additions and the evaluation of simple nonlinearities, which are in this case rectified linear units.
Both networks were deployed on a low-cost \SI{32}{\bit} microcontroller (SAMD21 Cortex-M0+) with \SI{32}{\kilo\byte} of RAM and \SI{48}{\mega\hertz} clock speed.
The evaluation time of the $\mathcal{N}_{6,6}$-controller for the oscillating masses example was \SI{1.6}{\milli\second} and the code required \SI{23.2}{\kilo\byte} of memory.
The code of the $\mathcal{N}_{10,6}$-controller for the inverted pendulum had a slightly larger memory footprint of \SI{24.6}{\kilo\byte} and the evaluation time was \SI{3.8}{\milli\second}.
The code was automatically generated with the open-source toolbox \emph{edgeAI} \cite{karg2018}.

\subsection{Binary Search Trees}

The use of Binary Search Trees can reduce the memory requirements and especially the evaluation time of the explicit MPC solutions~\cite{tondel2003evaluation} compared to a standard explicit MPC implementation.
However, we have not included the corresponding results for BSTs for the given examples since computing them for $N \geq 3$ with the toolbox MPT3~\cite{MPT3} was intractable.
For shorter horizons, the BST led to a memory footprint reduction around \SI{25}{\percent}.
For instance for the inverted pendulum on a cart with horizon $N=2$ the BST led to a reduction of \SI{21.9}{\percent}, but it was not possible to solve problems with longer horizons.

\section{Conclusions and future work}\label{sec:conclusions}

We have shown that explicit MPC formulations can be exactly represented by deep neural networks with rectifier units as activation functions and included explicit bounds on the dimensions of the required neural networks. The choice of deep networks is especially interesting for the representation of explicit MPC laws as the number of regions that deep networks can represent grows exponentially with their depth. 
This notion was exploited to propose an approximation method for explicit MPC solutions.

Stochastic verification techniques have been used to ensure constraint satisfaction if the neural network is used to approximate, and not to exactly represent, the explicit MPC law within a safe set.
Simulation results show that the proposed deep learning-based explicit MPC achieves better performance than other approximate explicit MPC methods with significantly smaller memory requirements.
This significant reduction of the memory footprint enabled the deployment of the proposed deep neural network controllers on a low-power embedded device with constrained resources.

Future work includes the design of stability guaranteeing formulations and the computation of safe sets which allow to predict the probability and the magnitude of constraint violations via gaussian process regression.
\ifCLASSOPTIONcaptionsoff
  \newpage
\fi



%
%
%

\bibliographystyle{IEEEtran}
\bibliography{refs_cdc}

%
\vspace{2cm}
\begin{IEEEbiography}[{\includegraphics[width=1in,height=1.25in,clip,keepaspectratio]{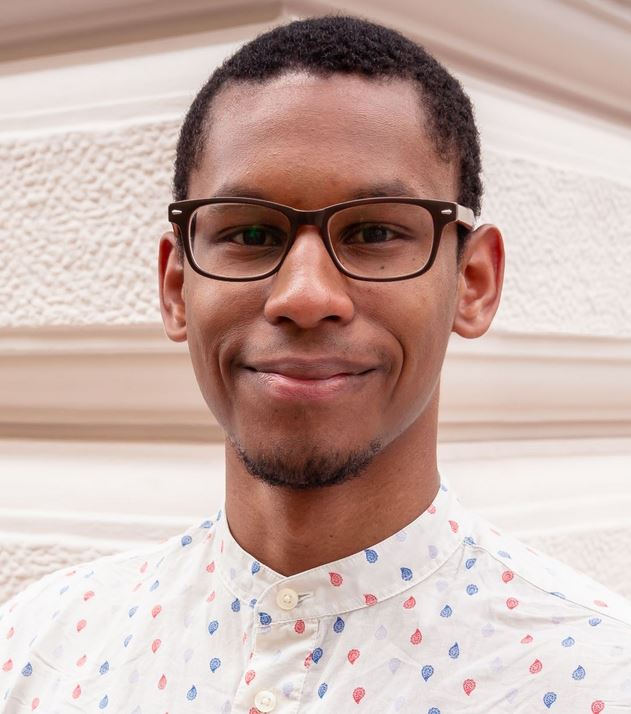}}]{Benjamin Karg}
was born in Burglengenfeld, Germany, in 1992.
He received the B.Eng. degree in mechanical engineering from Ostbayerische Technische Hochschule Regensburg, Regensburg, Bavaria, Germany, in 2015, and his M.Sc. degree in systems engineering and engineering cybernetics from Otto-von-Guericke Universität, Magdeburg, Saxony-Anhalt, Germany, in 2017.

He currently works as a research assistant at the laboratory "Internet of Things for Smart Buildings", Technische Universität Berlin, Germany, to pursue his PhD.
He is also member of the Einstein Center for Digital Future.

His research is focused on control engineering, artificial intelligence and edge computing for IoT-enabled cyber-physical systems.
\end{IEEEbiography}
\vspace{-2cm}
\begin{IEEEbiography}[{\includegraphics[width=1in,height=1.25in,clip,keepaspectratio]{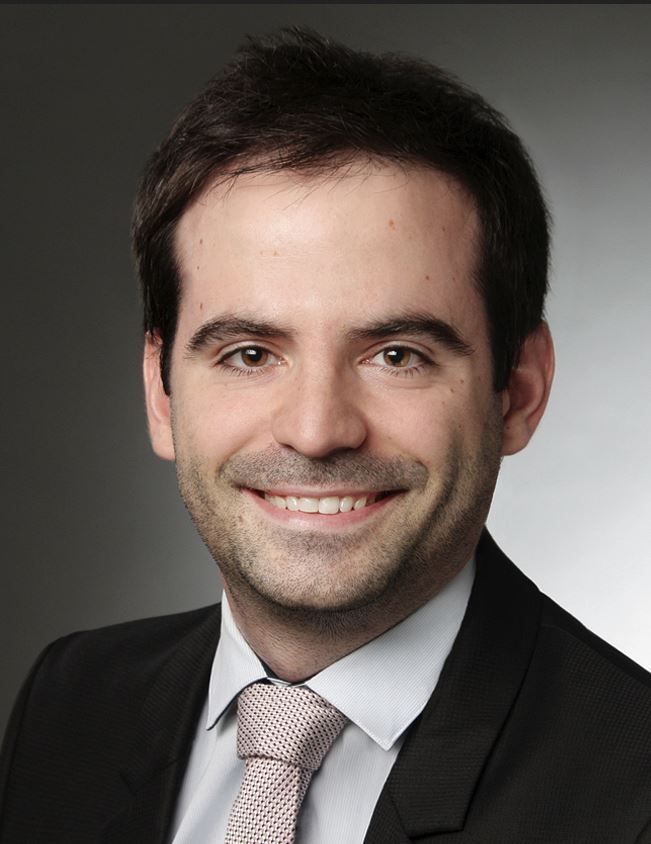}}]{Sergio Lucia}
(M’16) received the M.Sc. degree
in electrical engineering from the University of
Zaragoza, Zaragoza, Spain, in 2010, and the
Dr. Ing. degree in optimization and automatic
control from the Technical University of Dortmund,
Dortmund, Germany, in 2014.
He joined the Otto-von-Guericke Universität Magdeburg and visited the Massachusetts
Institute of Technology as a Postdoctoral Fellow.

Since May 2017, he has been an Assistant Professor
and Chair with the Laboratory of “Internet
of Things for Smart Buildings”, Technische Universität Berlin, Berlin,
Germany, and with Einstein Center Digital Future, Berlin.
His research interests include decision-making under uncertainty, distributed control, as well as the interplay between machine learning techniques and control theory.

Dr. Lucia is currently Associate Editor of the Journal of Process Control.
\end{IEEEbiography}






\end{document}